\documentclass[12pt,reqno,twoside]{amsart}

\usepackage{amsfonts}
\usepackage{amsmath,amssymb,amsthm,amsthm}
\usepackage{mathtools}
\usepackage{dsfont}
\usepackage{etoolbox}
\usepackage{color}
\usepackage{bm}

\usepackage[dvips,bottom=1.4in,right=1in,top=1in, left=1in]{geometry}

\makeatletter
\def\eqnarray{\stepcounter{equation}\let\@currentlabel=\theequation
\global\@eqnswtrue
\tabskip\@centering\let\\=\@eqncr
$$\halign to \displaywidth\bgroup\hfil\global\@eqcnt\z@
  $\displaystyle\tabskip\z@{##}$&\global\@eqcnt\@ne
  \hfil$\displaystyle{{}##{}}$\hfil
  &\global\@eqcnt\tw@ $\displaystyle{##}$\hfil
  \tabskip\@centering&\llap{##}\tabskip\z@\cr}

\def\endeqnarray{\@@eqncr\egroup
      \global\advance\c@equation\m@ne$$\global\@ignoretrue}

\numberwithin{equation}{section}

\def\ms{\medskip}

\def\ub{\underbar}
\def\noi{\noindent}

\def\khat{{\hat k}}

\def\utilde{{\tilde u}}
\def\vtilde{{\tilde v}}

\def\fvec{{\bf f}}

\def\kvec{{\bf k}}

\def\nvec{{\bf n}}

\def\rvec{{\bf r}}

\def\uvec{{\bf u}}

\def\xvec{{\bf x}}

\def\Cmat{{\bf C}}

\def\Kmat{{\bf K}}
\def\Lmat{{\bf L}}
\def\Mmat{{\bf M}}

\def\Grd{\nabla}
\def\Div{\nabla \cdot}

\def\pmb#1{\setbox0=\hbox{$#1$}%
             \kern-.027em\copy0\kern-\wd0
             \kern+.009em\copy0\kern-\wd0
             \kern+.009em\copy0\kern-\wd0
             \kern+.009em\copy0\kern-\wd0
             \kern+.009em\copy0\kern-\wd0
             \kern+.009em\copy0\kern-\wd0
             \kern+.009em\copy0\kern-\wd0
             \kern-.045em\raise+.012em\copy0\kern-\wd0
             \kern+.009em\raise+.012em\copy0\kern-\wd0
             \kern+.009em\raise+.012em\copy0\kern-\wd0
             \kern+.009em\raise-.012em\copy0\kern-\wd0
             \kern+.009em\raise-.012em\copy0\kern-\wd0
             \kern-.018em\copy0\kern-\wd0\raise-.012em\box0}
\def\Pmb#1{\setbox0=\hbox{$#1$}%
             \kern-.033em\copy0\kern-\wd0
             \kern+.011em\copy0\kern-\wd0
             \kern+.011em\copy0\kern-\wd0
             \kern+.011em\copy0\kern-\wd0
             \kern+.011em\copy0\kern-\wd0
             \kern+.011em\copy0\kern-\wd0
             \kern+.011em\copy0\kern-\wd0
             \kern-.055em\raise+.015em\copy0\kern-\wd0
             \kern+.011em\raise+.015em\copy0\kern-\wd0
             \kern+.011em\raise+.015em\copy0\kern-\wd0
             \kern+.011em\raise-.015em\copy0\kern-\wd0
             \kern+.011em\raise-.015em\copy0\kern-\wd0
             \kern-.022em\copy0\kern-\wd0\raise-.015em\box0}
\vfill

\usepackage[textsize=small]{todonotes}
\title[]{Revisiting Calderon's Problem}

\author{Rainald L\"ohner}
\address{Center for Computational Fluid Dynamics, College of Science, 
George Mason University, Fairfax, VA 22030-4444, USA.}
\email{rlohner@gmu.edu}

\author{Harbir Antil}
\address{Department of Mathematical Sciences, 
George Mason University, Fairfax, VA 22030, USA.}
\email{hantil@gmu.edu}

\thanks{The work of the second author is partially supported by NSF 
grants DMS-1818772 and DMS-1913004 and the Air Force Office of 
Scientific Research under Award NO: FA9550-19-1-0036.}

\keywords{Calderon Problem, Inverse Problems, Parameter Estimation, Optimization Problem,
  Heat Conduction, Finite Element Method.}

\subjclass[2010]{
  35R30,  	
  49N45,  	
  65N21,  	
  49K20,  	
  49J20,  	
  65N30.  	
}

\begin{document}

\begin{abstract}
A finite element code for heat conduction, together with an
adjoint solver and a suite of optimization tools was applied
for the solution of Calderon's problem. One of the questions 
whose answer was sought was whether the solution to these 
problems is unique and obtainable. The results to date show 
that while the optimization procedure is able to obtain spatial
distributions of the conductivity $k$ that reduce the cost function
significantly, the resulting conductivity $k$ is still significantly 
different from the target distribution sought. While the normal fluxes 
recovered are very close to the prescribed ones, the tangential fluxes 
can differ considerably.
\end{abstract}

\maketitle

\section{Introduction}

The problem of trying to determine the material properties of
a domain from boundary information is common to many fields.
To mention just a few:
mining (e.g. prospecting for oil and gas), medicine (e.g. trying
to infer tissue properties), and engineering (e.g. trying to determine
the existence and location of fissures). 

From an abstract setting, it would seem that this is an ill-posed
problem. After all, if we think of atoms, granules or some polygonal
(e.g. finite element [FEM]) subdivision of space, the amount of data given
resides in a space of one dimension less than the data sought. If we
think of cuboid domain in $d$ dimensions with $N^d$ subdivisions, 
the amount of information/ data given is of $O(N^{d-1})$ while the 
data sought is of $O(N^d)$.

Another aspect that would seem to indicate that this is an ill-posed
problem is the possibility that many possible spatial distributions
of material properties could yield very similar or equal boundary
values. That this is indeed the case for some problems is shown below.

On the other hand, the propagation of physical properties (e.g.
temperature, displacements, electrical currents, etc.) through the
domain obeys physical conservation laws, i.e. some partial differential
equations (PDEs). This implies that the material properties that 
can give rise to the data measured on the boundary are restricted by 
these conservation laws, i.e. are bounded. This would indicate that 
perhaps -~due to these restrictions~- the problem is not as ill-posed 
a initially thought.

\section{Calderon's Problem}

In the following, we will consider conservation laws of the form:

\begin{equation}\label{eq:conslaw}
	\Div \fvec = 0  \, , 
\end{equation}
with
\begin{equation}\label{eq:law} 
	\fvec = k \Grd u  \, ,
\end{equation}	
implying
\begin{equation}\label{eq:Diffeq} 
	\Div k \Grd u = 0     \, , \quad \mbox{in} ~~\Omega \, , 
\end{equation}	
where $u, k, \fvec$ denote the unknowns, material property and
flux, and $\Omega \subset \mathbb{R}^d$, with $d \ge 1$ the domain 
considered with boundary $\Gamma$. 
For heat transfer problems, $u$ is the temperature, $k$ the
conductivity, $\fvec$ the heat flux, and \eqref{eq:law} Fourier's law.
For electrical currents, $u$ is the voltage, $k$ the
resistivity, $\fvec$ the electrical current, and \eqref{eq:law} Ohm's law.

Eqn.~\eqref{eq:Diffeq} needs appropriate boundary conditions in order to be 
uniquely solvable. On the boundary $\Gamma$ we can prescribe 
$u$ (Dirichlet), $f_n=\nvec \cdot \fvec$ (Neumann) or a combination 
of both (Robin), i.e.
\begin{equation}
	\alpha u + \beta f_n = b \, , \quad \mbox{on} ~~\Gamma \, ,
\end{equation}	

Calderon's problem \cite{APCalderon_1980a} may then be formulated as follows:
\begin{enumerate}
	\item[-] Given a PDE of the form of \eqref{eq:Diffeq};
	\item[-] Given {\bf both} $u$ and $f_n$ on $\Gamma$;
	\item[-] Determine $k(\xvec)$ in $\Omega$.
\end{enumerate}
Ever since Calderon first formulated this problem in 1980 
\cite{APCalderon_1980a},
several proofs of existence, uniqueness and solvability have been
given \cite{GAlessandrini_1988a,GAlessandrini_1990a,RMBrown_GAUhlmann_1997a,
JFeldman_MSalo_GUhlmann_2019a}. 
However, as will be seen, the problem remains difficult. 
It is also interesting to observe that there is apparently no
`canonical test problem' or `canonical test suite of problems' for
this class of problems. This is in sharp contrast to other engineering
disciplines such as aerodynamics, hydrodynamics and electromagnetics,
where standard test problems have evolved over the years.

\section{The 1-D Case}

In 1-D, \eqref{eq:conslaw}-\eqref{eq:Diffeq} reduce to:
\begin{equation}\label{eq:1Dconslaw}
	f=f_c \, , \quad k {du \over dx} = f_c \, , 
\end{equation}	
where $f_c$ is the constant heat flux. This implies that only one
measurement point is required for $f_n$ (note also: in 1-D $f_n = f$). 
Given a domain $0 \le x \le L$ and $k > 0$, we can integrate 
\begin{equation}\label{eq:1Dsys}
	{du \over dx} = { f_c \over k } \, , 
\end{equation}
implying 
\begin{equation}\label{eq:1Dints}
	u(L) - u(0) = f_c \int_0^L { 1 \over k } dx \, ,
\end{equation}	
or:
\begin{equation}\label{eq:1Dksolve}
	\int_0^L { 1 \over k } dx = {{u(L) - u(0)} \over f_c}  \, . 
\end{equation}                                               
Calderon's problem in 1-D then reduces to determining $k$ given
$u(0), u(L)$, and $f_c$. This is unique for constant $k(x)=k_c$, but not
for arbitrary distributions of $k(x)$. If we consider $n$ 
regions $\Delta x_i, \, i=1, \dots, n$ with different $k_i$, then 
all that is required is:
\begin{equation}\label{eq:1DksolveD}
	\sum_{i=1}^n {{\Delta x_i} \over k_i} = {{u(L) - u(0)} \over f_c} \, ,
\end{equation}   
which clearly allows for infinitely many possible solutions.
To illustrate this, consider the simple case with:
$L=1, u(0)=1, u(1)=0, f_c=-1$. Assuming 4 equal regions of 
$\Delta x=0.25$, \eqref{eq:1DksolveD} implies:
\begin{equation}\label{eq:1Dex}
	 \sum_{i=1}^4 {1 \over k_i} = 4 \, , 
\end{equation}

\begin{figure}[h!]
	\centering
	\includegraphics[width=0.45\textwidth]{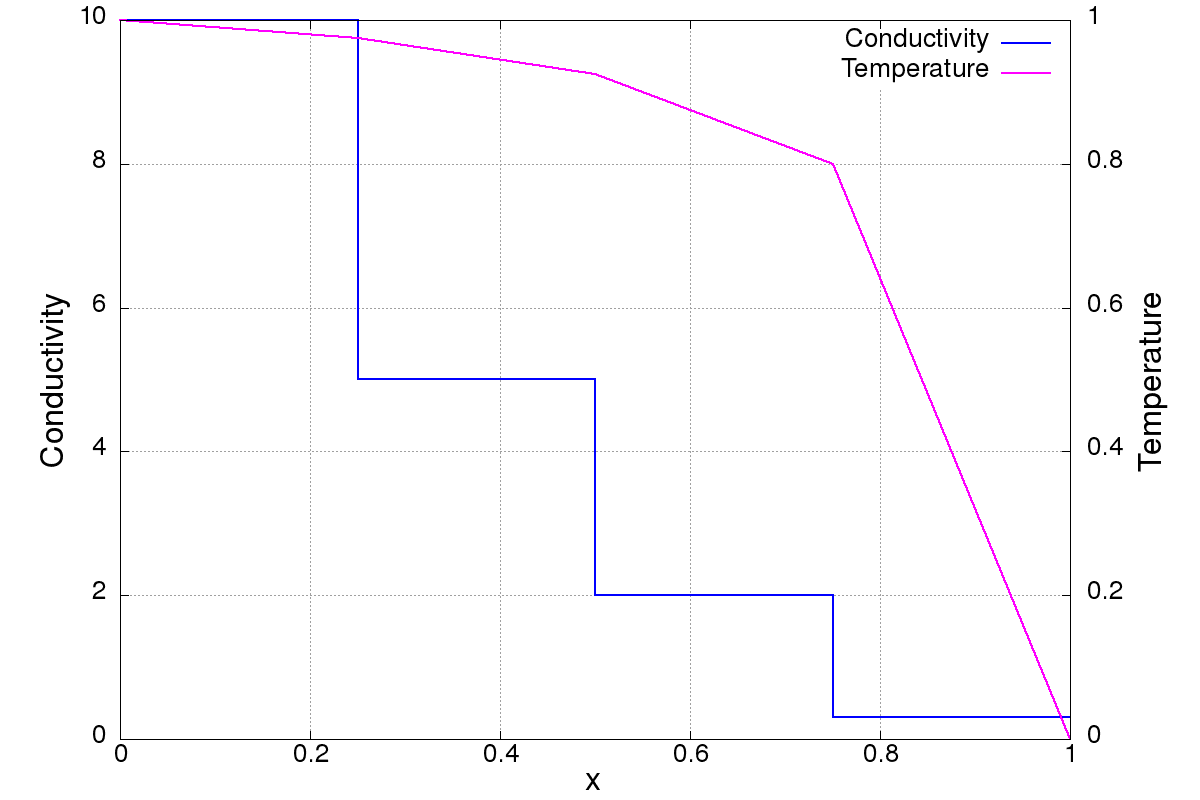}
	\includegraphics[width=0.45\textwidth]{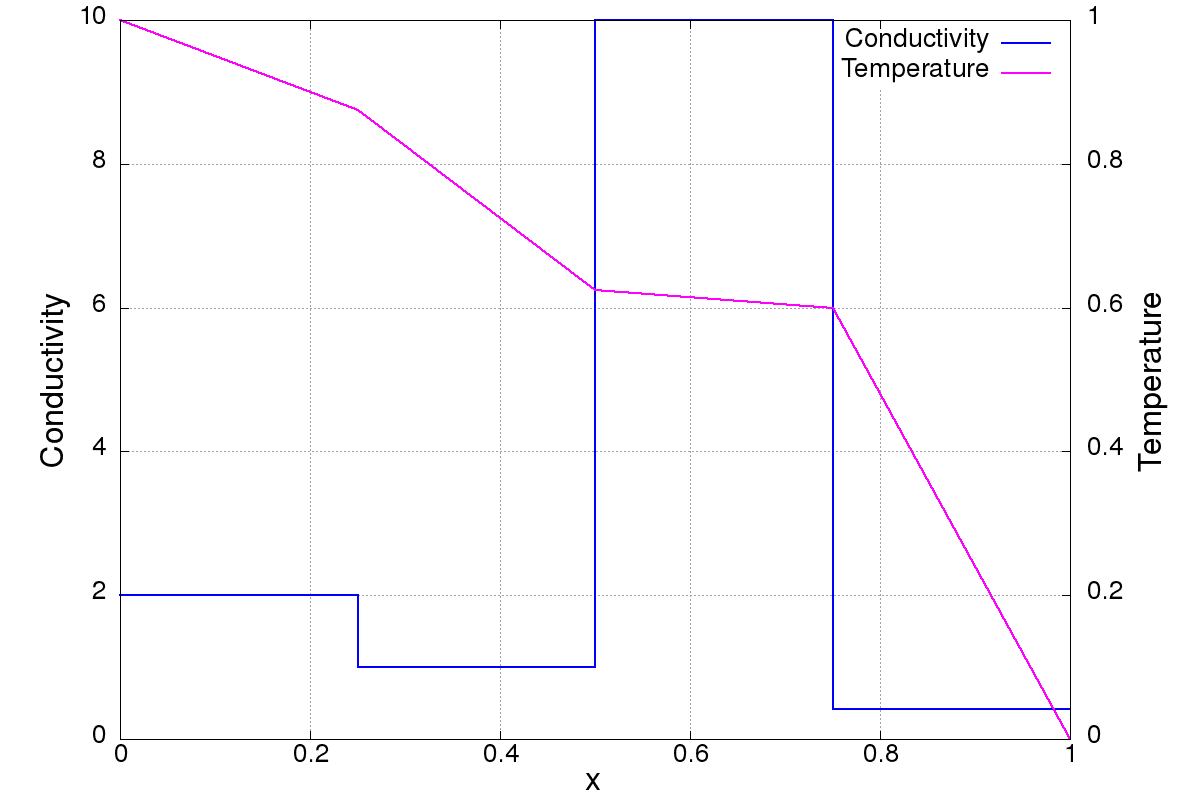}
	\caption{\label{f:1D} 1D Case (a), (b): Possible Conductivities and Resulting Temperatures 
	for $f_c=-1$.}
\end{figure}

\begin{figure}[h!]
	\centering
	\includegraphics[width=0.45\textwidth]{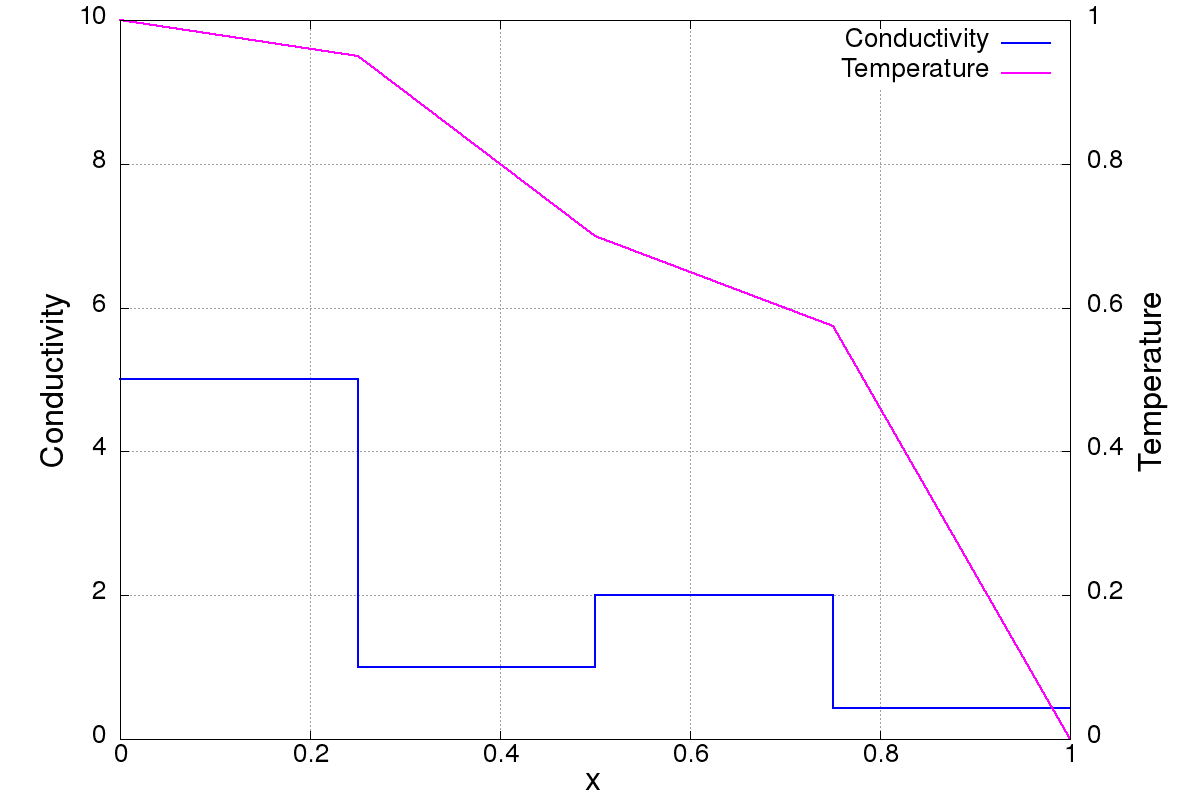}
	\includegraphics[width=0.45\textwidth]{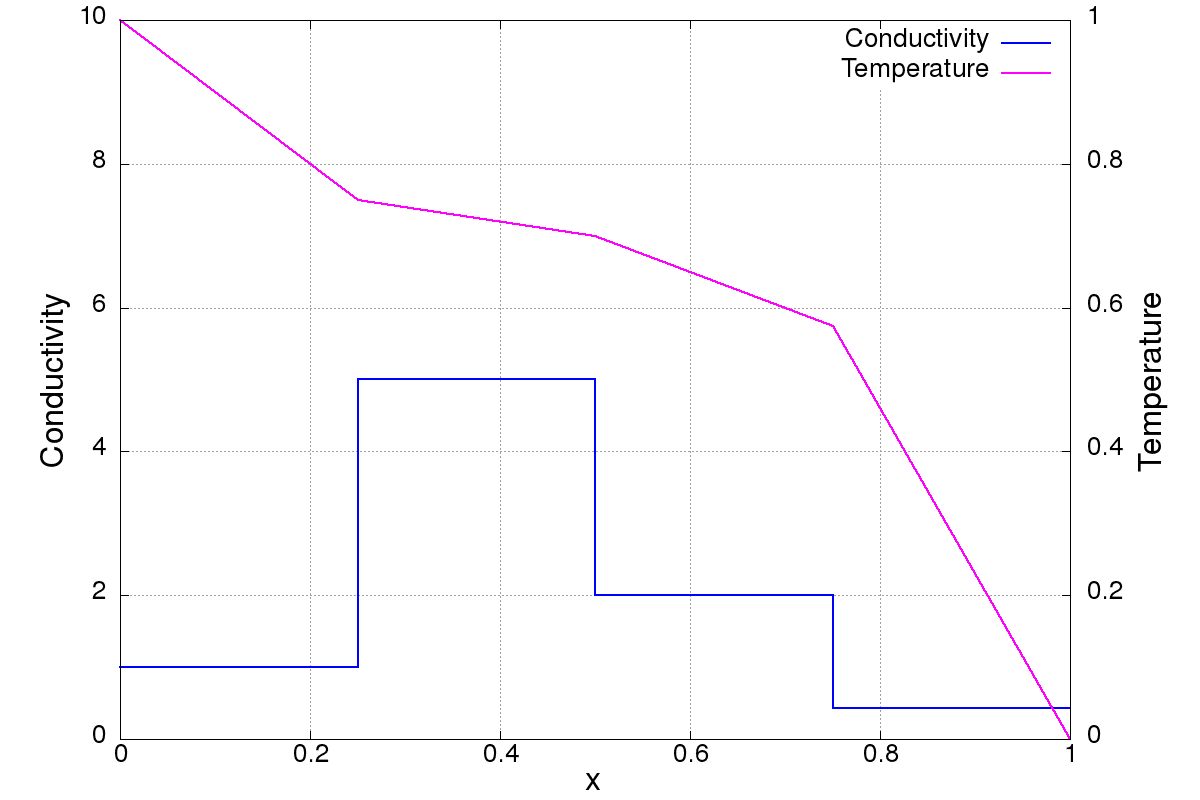}
	\caption{\label{f:1D2} 1D Case (c), (d): Possible Conductivities and Resulting Temperatures 
	for $f_c=-1$.}
\end{figure}

\noi
which may be realized, e.g. via
\begin{enumerate}
	\item[a)] $k_1 = k_2 = k_3 = k_4 = 1$;
	\item[b)] $k_1 = k_2 = 2; k_3 = k_4 = 2/3 $;
	\item[c)] $k_1 = 10, k_2 = 5, k_3 = 2, k_4 = 10/32 $;
	\item[d)] ...
\end{enumerate}
Moreover, the spatial sequence of $k$'s can be changed without
affecting the resulting flux, so that even highly oscillatory
distributions of $k$ are admissible. Four of these possibilities
have been plotted in Figure~\ref{f:1D} (a), (b) and Figure~\ref{f:1D2} (c), (d).    

Why would Calderon's problem then be easier to solve or perhaps uniquely
solvable in higher dimensions? Consider again the 1-D case, and
assume that this is a bar/hexahedral domain with equal material 
properties along the $y,z$ directions, as well as Neumann boundary 
conditions at the $y,z$ boundaries. In the 1-D case, we have 
restricted the data input to only 3 items: a constant temperature 
for planes $x=0, x=L$, and the heat flux at $x=0$. In the 3-D case,
we could measure the temperature and fluxes at all boundaries
(and in particular the $y,z$ boundaries along the $x$ axis, and we
could impose a temperature that is not constant in $y,z$ for the
planes $x=0, x=L$. This would yield much more information than in the
1-D case, and may lead to a solvable problem.

\begin{figure}[h!]
	\centering
	\includegraphics[width=0.7\textwidth]{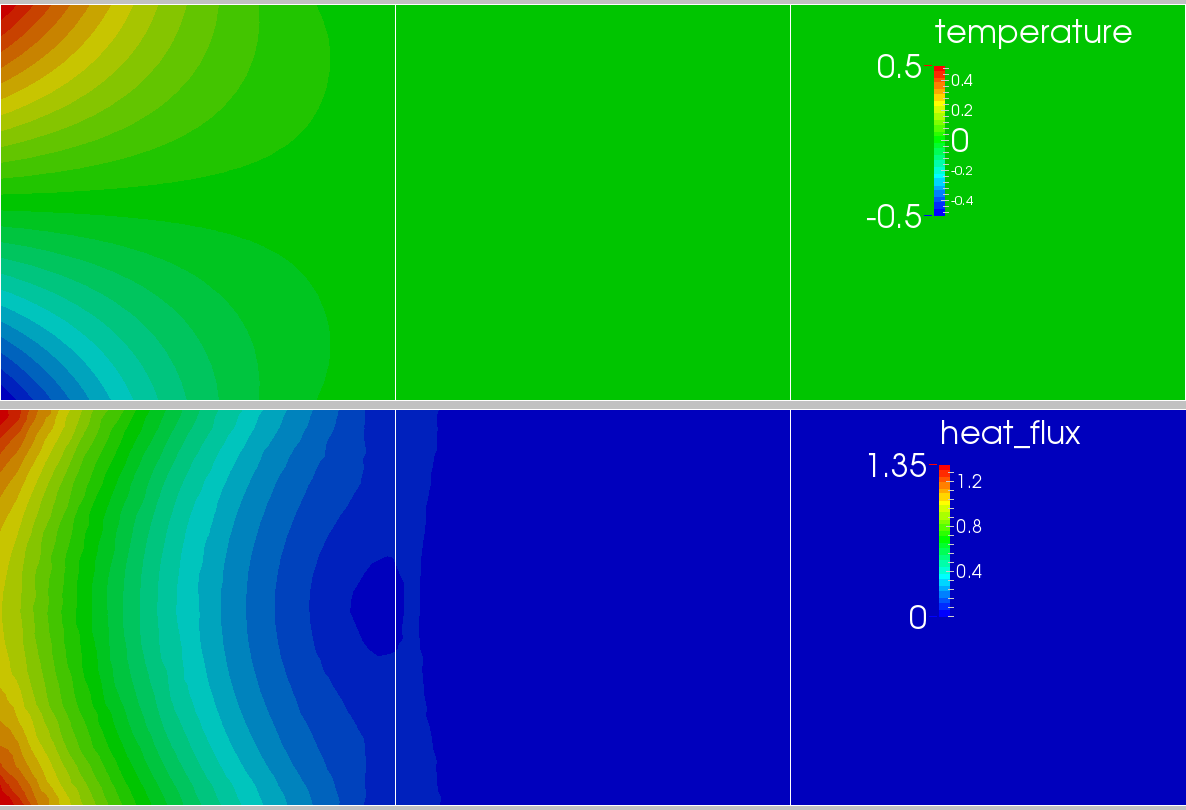}
	\caption{\label{f:2D} 2D Case With Small Values for $\nabla u$.}
\end{figure}

On the other hand, one may encounter difficulties for the 
multidimensional cases when regions with $\nabla u \approx 0$ 
are present. In these cases it is irrelevant which value of 
$k$ is used, as the flux is vanishingly small anyhow. As an example, 
consider the 2-D case $0 \le x \le 3$, $0 \le y \le 1$ shown in 
Figure~\ref{f:2D}. The domain is divided into 3 regions and the 
the conductivity in these three regions 
is set to $k_{left}=1.0, k_{middle}=10.1, k_{right}=1.0$. The value of
$u$ on the boundary was set to $u|_{\Gamma}=0$, except for
$x \le 1$ where $u|_{\Gamma}=(y-0.5)(x-1)^2$. As one can see from the
temperature and the fluxes, in the right part of the domain very small
gradients of $u$ are obtained. This would indicate that one has to be
careful when `conducting a measurement', imposing values of $u, f_n$ that
yield useful information about the spatial distribution of $k$.         

\section{Calderon's Problem and Optimization}

Calderon's problem may be formulated as an optimization problem 
with PDE constraints for $k(\xvec)$ as follows: 
Given $m$ measurements with data 
$\{ (u_m, f_m) \}_{m}$ on the boundary, minimize the cost functional:
\begin{equation}\label{eq:cost}
	I(u_m,k) = 
		\frac{1}{2} \sum_m \int_{\Gamma} ( f_m - k ~\nvec \cdot \nabla u_m )^2 d \Gamma \, , 
\end{equation}                                                   
subject to PDE constraints
\begin{equation}\label{eq:PDE}
	R(u_m,k) = \Div k \Grd u_m = 0     \, , \quad \mbox{in } \Omega  
\end{equation}	
with boundary condition
\begin{equation}\label{eq:BC}
	u_m = u_m^0 \, , \quad \mbox{on } \Gamma \, . 
\end{equation}
In the most general case, one may consider a spatial subdivision 
of space (e.g. a FEM discretization) and functions with local support 
so that:
\begin{equation*}
	k(\xvec) = \sum_i N^i(\xvec) \khat_i \, , 	
\end{equation*}	
where $N^i(\xvec)$ are known functions (e.g. FEM shape functions) and
$\khat_i$ the discrete values (i.e. degrees of freedom, design
parameters) sought.
At the core of any optimization tool is the cost function evaluation.
In the present case, this implies invoking a field solver $m$ times
for \eqref{eq:PDE}, in the weak form, followed by boundary evaluation 
for \eqref{eq:cost}.
As we are seeking domain information, the number of design 
variables may be very large. This implies that non-gradient based
procedures (e.g. genetic algorithms), or finite diffence based
gradient techniques may be prohibitively expensive.

\section{Optimization via Adjoints}

The purpose of this section is to use the formal Lagrangian approach 
to derive the necessary optimality conditions for 
\eqref{eq:cost}-\eqref{eq:BC}. The cost functional can be extended to 
the Lagrangian functional
\begin{equation}
L(u_m,k,\utilde_m,\tilde{v}_m) = I(u_m,k) 
- \sum_m \int_{\Omega} \utilde_m  \Div k \Grd u_m d\Omega 
+ \sum_m \int_{\Gamma} \vtilde_m ( u_m - u^0_m) d\Gamma \, , 
\end{equation}
where $\utilde_m, \vtilde_m$ are the Lagrange multipliers (adjoints) 
corresponding to the state equation \eqref{eq:PDE} and the boundary 
condition \eqref{eq:BC}, respectively. 
Repeated integration by parts leads to:
\begin{align}\label{eq:LagIP}
\begin{aligned}
L(u_m,k,\utilde_m,\tilde{v}_m) &= I(u_m,k) 
+ \sum_m \int_{\Omega} k \nabla \utilde_m \Grd u_m d\Omega \\
&\quad- \sum_m \int_{\Gamma} \utilde_m k ~\nvec \cdot \nabla u_m d\Gamma
+ \sum_m \int_{\Gamma} \vtilde_m ( u_m - u^0_m) d\Gamma  \, ,
\end{aligned}
\end{align}
\begin{align*}
\begin{aligned}
L(u_m,k,\utilde_m,\tilde{v}_m) = I(u_m,k) 
- \sum_m \int_{\Omega} u_m \Div k \Grd \utilde_m d\Omega 
+ \sum_m \int_{\Gamma} u_m k ~\nvec \cdot \nabla \utilde_m d\Gamma \\
- \sum_m \int_{\Gamma} \utilde_m k ~\nvec \cdot \nabla u_m d\Gamma
+ \sum_m \int_{\Gamma} \vtilde_m ( u_m - u^0_m) d\Gamma \, . 
\end{aligned}
\end{align*}

\noi
Considering a variation with respect to $u_m$ in a direction $h$ results in:
\begin{align*}
\begin{aligned}
D_{u_m} L(u_m,k,\utilde_m,\tilde{v}_m) h = 
- \int_{\Gamma} \left[ ( f_m - k ~\nvec \cdot \nabla u_m ) 
                        k ~\nvec \cdot \nabla h 
- h k ~\nvec \cdot \nabla \utilde_m \right. \\
\quad \left.
+ \utilde_m k ~\nvec \cdot \nabla h - h \tilde{v}_m  \right] d \Gamma
- \int_{\Omega} h \Div k \Grd \utilde_m d\Omega  = 0 \, .
\end{aligned}
\end{align*}

Selecting $h$ so that $h$ and all its derivatives are zero on 
$\Gamma$ yields the \emph{adjoint equation} in the domain $\Omega$:
\begin{equation}\label{eq:adjPDE}
	\Div k \Grd \utilde_m = 0  \, . 
\end{equation}	
Choosing $h$ such that $h = 0$ on $\Gamma$ but 
$ k ~\nvec \cdot \nabla h \neq 0$ yields the boundary conditions 
for $\tilde{u}_m$:
\begin{equation}\label{eq:adjBC}
	\utilde_m = - \left( f_m - k ~\nvec \cdot \nabla u_m \right)
                                                \, . 
\end{equation}       
Notice that by letting $h\neq 0$ on $\Gamma$ gives a compatibility condition 
between $\tilde{v}_m$ and $\tilde{u}_m$.
Finally, performing a variation of $L$, in \eqref{eq:LagIP}, 
with respect to $k$ in a direction $h$ leads to:
\begin{align}
\begin{aligned}
	D_k L(u_m,k,\utilde_m,\tilde{v}_m) h  =
\sum_m \left[ \int_\Gamma - \left(f_m - k ~\nvec \cdot \nabla u_m \right) 
                                 \nvec \cdot \nabla u_m h d\Gamma
+ \int_{\Omega} h \nabla u_m \cdot \nabla \utilde_m d\Omega 
\right.  \\
 \left.
-        \int_{\Gamma} \utilde_m \nvec \cdot \nabla u_m h d\Gamma
                                     \right] \, , 
\end{aligned}
\end{align}
which, after using \eqref{eq:adjBC}, is equivalent to:
\begin{equation}\label{eq:kVar}
	D_k L(u_m,k,\utilde_m,\tilde{v}_m) h  = \sum_m \left[
              \int_{\Omega} \nabla u_m \cdot \nabla \utilde_m h d\Omega 
                                     \right] \, .
\end{equation}                                     
The consequences of this rearrangement are profound:
\begin{enumerate}
\item[-] The gradient of $L$, $I$ (cf.~\ref{eq:kVar}) with respect 
to $k$ may be obtained by solving $m$ forward \eqref{eq:PDE}-\eqref{eq:BC} 
and adjoint \eqref{eq:adjPDE}-\eqref{eq:adjBC} problems; i.e.
\item[-] The cost for the evaluation of gradients in \eqref{eq:kVar} is 
{\bf independent of the number of variables used for $k$} (!). This implies
that the material coefficients $k$ sought may be described by a very large
parameter set (for example a constant value for each Finite Element),
something that would be computationally prohibitively expensive for 
methods that obtain gradients via finite differences of the objective
function.
\end{enumerate}

\noi
A gradient descent based optimization cycle using the adjoint approach 
is then composed of the following steps:
\begin{enumerate}
\item[-] For each measurement $m$:
	\begin{enumerate}
	\item[-] With current $k$: solve \eqref{eq:PDE}-\eqref{eq:BC} for 
                 the field variable $\rightarrow u_m$
	\item[-] With current $k$ and $u_m$: solve 
                 \eqref{eq:adjPDE}-\eqref{eq:adjBC} for the adjoint field
		 variable $\rightarrow \utilde_m$
	\item[-] With $u_m, \utilde_m$: obtain gradients for each $m$
		 $\rightarrow I(u_m)_{,k}$
	\end{enumerate}
\item[-] Sum up the gradients $\rightarrow I_{,k} = \sum_m I(u_m)_{,k}$
\item[-] Smooth gradients (regularization) $\rightarrow I^s_{,k}$
\item[-] Update $k_{new} = k_{old} - \alpha I^s_{,k}$, where $\alpha$ 
         is the step-length.
\end{enumerate}
Here and subsequently $I_{,k}:={\partial I}/{\partial k}$.

\section{Smoothing of Gradients}

The gradients of the cost function with respect to $k$ allow for 
oscillatory solutions. One must therefore smooth or `regularize' the
spatial distribution. This happens naturally when using few degrees of
freedom, i.e. when $k$ is defined via other spatial shape functions
(e.g. larger spatial regions of piecewise constant $k$). As the 
(possibly oscillatory) gradients obtained in the (many) finite elements 
are averaged over spatial regions, an intrinsic smoothing occurs.
This is not the case if $k$ and the gradient are defined and evaluated
in each element separately, allowing for the largest degrees of
freedom in a mesh and hence the most accurate representation. 
Three different types of smoothing or `regularization' were considered.
All of them start by performing a volume averaging from elements 
to points: 
\begin{equation}\label{eq:kavpatch} 
 \khat_i = {{ \sum_{el} \khat_{el} V_{el} } \over
              { \sum_{el} V_{el} }} \, ,
\end{equation}              
where $\khat_i, \khat_{el}, V_{el} $ denote the value of $k$ at 
point $i$, as well as the values of $k$ in element $el$ and the 
volume of element $el$,  and the sum over all the elements
surrounding point $i$.

\subsection{Simple Point/Element/Point Averaging}

In this case, the values of $k$ are cycled between elements and points.
When going from point values to element values, a simple average is
taken:
\begin{equation}\label{eq:SPEA}
	\khat_{el} = { 1 \over n_{el} } \sum_i \khat_i \, ,
\end{equation}	
where $n_{el}$ denotes the number of nodes (degrees of freedom) of
an element and the sum extends over all the nodes of the element.
After obtaining the new element values via \eqref{eq:SPEA} the point 
averages are again evaluated via \eqref{eq:kavpatch}. 
This form of averaging is very crude, but works surprisingly well.

\subsection{$H^1$ (Weak) Laplacian Smoothing}
In this case, the initial values $k_0$ obtained for $k$ are smoothed via:

\begin{equation}\label{eq:H1smooth}
	\left[ 1 - \lambda_l \nabla^2 \right] k = k_0 \, . 
\end{equation}	
Here $\lambda$ is a free parameter which may be problem and mesh 
dependent (its dimensional value is length squared). Discretization 
(of the weak form of \eqref{eq:H1smooth}) via FEM yields:
\begin{equation}\label{eq:H1smoothD}
 \left[ \Mmat_c + \lambda_l \Kmat \right] \kvec = 
          \Mmat_c \kvec_0 \, ,
\end{equation}          
where $\Mmat_c, \Kmat$ denote the consistent mass matrix and the 
stiffness matrix obtained for the Laplacian operator.

\subsection{$H^1$ (Weak) Pseudo-Laplacian Smoothing}
One can avoid the dimensional dependency of $\lambda$ in 
\eqref{eq:H1smooth} by smoothing via:
\begin{equation}\label{eq:PH1smooth}
	\left[ 1 - \lambda_{pl} \nabla \cdot h_{fem}^2 \nabla \right] k = k_0 \, ,
\end{equation}	
where $h_{fem}$ is a characteristic element size (e.g. average of element
sides, average of element normals, or any other sensible measure of
element size). For linear elements, one can show that this is equivalent to:
\begin{equation}\label{eq:PH1smoothD}
	\left[ \Mmat_c + \lambda_{pl} \left(\Mmat_l - \Mmat_c \right) \right] 
	      \kvec = \Mmat_c \kvec_0 \, ,
\end{equation}	      
where $\Mmat_l$ denotes the lumped mass matrix \cite{RLoehner_2008a}. 
In the examples shown below this form of smoothing was used for the 
gradients, setting $\lambda_{pl}=0.05$.

\subsection{Relaxation}

As the main aim of the smoothing given by 
\eqref{eq:H1smooth}-\eqref{eq:PH1smoothD} 
is the removal of high frequency modes, simple relaxation procedures that
do not require the inversion of a large matrix offer an effective
way to proceed. Eqns.~\eqref{eq:H1smoothD}, 
\eqref{eq:PH1smoothD} may be recast as:
\begin{equation*}
	\Lmat \uvec = \rvec \, . 
\end{equation*}	

\noi
This may be interpreted as the steady result of the following transient
problem:
\begin{equation*}
	\Cmat \uvec_{,\tau} + \Lmat \uvec = \rvec \, . 
\end{equation*}	
where $\tau$ denotes a pseudo-time and $\Cmat$ contains 
the diagonal entries of $\Lmat$. This is solved via explicit integration:
\begin{equation*}
\Cmat \left( \uvec^{n+1} - \uvec^n \right) =
\delta \tau \left(  \rvec - \Lmat \uvec^n \right) \,  .
\end{equation*}
A value of $\delta\tau=0.8$ was set for all numerical examples in the next 
section. It was observed that the number of explicit steps required for 
effective smoothing is of the order of 10. For 3-D cases, it is 
far more economical to use such a relaxation procedure than inverting 
the matrix $\Lmat$.

\subsection{Gradient Projection}

Intrigued by some of the results obtained, and in order to check whether
the gradients obtained via the adjoint method were correctly computed, 
the option of projecting the gradients to a region discretization of 
$k$ based on volumes $V^{reg}_i$ with constant $k$ was implemented. 
Given the discretization:
 \begin{equation*}
	I_{,k}(\xvec) \approx \sum_i H_i \left( \hat{I^{reg}}_{,k} \right)_i \, ,
\end{equation*}	
where 
\[
	H_i(\xvec)= \left\{ \begin{array}{ll}
					1 & \forall \xvec \in    V^{reg}_i \, , \\
					0 & \forall \xvec \notin V^{reg}_i \, .
				\end{array}	
			  \right.
\]	   
and $(\hat{I^{reg}}_{,k})_i$ denote the (constant) values of the gradient
of $I$ with respect to $k$ in region $V^{reg}_i$, and given
$I_{,k}(\xvec)$ from the (fine) FEM mesh, the gradient is averaged over 
each volume as:
\[
	\left( \hat{I^{reg}}_{,k} \right)_i = 
           {{\int_{V^{reg}_i} I_{,k}(\xvec) dV} \over {\int_{V^{reg}_i} dV}} \, .        
\]
This gradient is then passed back to gradient-based optimizers to continue
the iteration towards the minimum of the cost function.

\section{Examples}

All the numerical examples are carried out using FEHEAT, a 
finite element code based on simple linear (beam), triangular
(plate) and tetrahedral (volume) elements with constant conductivity
per element. The optimization loops are steered via a
user-defined cost function and FEOPT, a general optimization code, for
the finite difference gradient-based optimization, or
via a simple shell-script for the adjoint-based optimization.
In all cases, a `target' distribution of $k(\xvec)$ and $u_{\Gamma}$
is generated with user-defined subroutines. The forward problem
is then solved, i.e. $u(\xvec)$ and $\fvec(\xvec)$ are obtained. This then
yields the `measurement data' $u_{\Gamma}, f_{\Gamma}$ that is 
used for the Calderon's problem. A series of tests are conducted to 
see if the conductivities obtained by the optimization algorithm
are dependent on the initial spatial distribution of the 
conductivity. The results showed no influence (all cases converged
to the same solution). Therefore, all optimization cases are started
from a constant $k(\xvec)=k_0$, and the convergence to the
`target' distribution of $k(\xvec)$, by minimizing the error with
respect to the `measured flux' given by \eqref{eq:cost}, is observed.
As stated in the introduction, it is possible that many different
spatial distributions of $k(\xvec)$ could minimize the error
measure given by \eqref{eq:cost}.

As it appears that no `canonical test problem' or `canonical 
test suite of problems' for this class of problems exists, 
an attempt has been made to start from simple cases and then proceed to
more difficult ones.

\subsection{Square}

\begin{figure}[h!]
	\centering
	\includegraphics[width=0.52\textwidth]{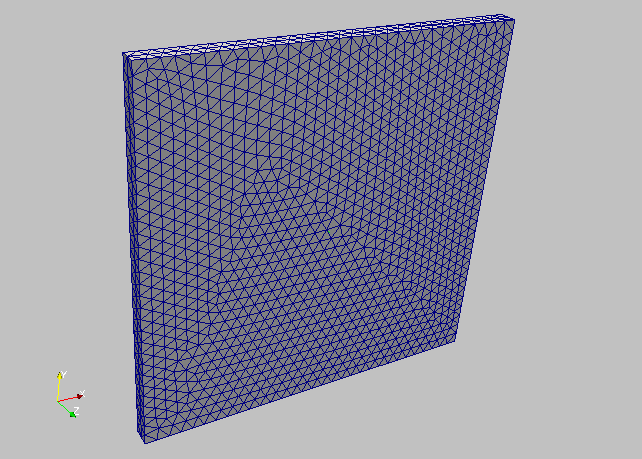}
	\caption{\label{f:squaremsh}Square: FEM Mesh Used.}
\end{figure}

This case is a 2-D case, but was run in 3-D. The domain and spatial
discretization are shown in Figure~\ref{f:squaremsh}. The dimensions are
$0 \le x \le 1, 0 \le y \le 1, 0 \le z \le 0.05$. This case was
run repeatedly, increasing the level of complexity and measurements.
The base mesh had elements of size approximately $h=0.025$ (the
ideal tetrahedron is not a space-filling element), which led to
12,461 elements and 3,347 points. The cases were also run on coarser
meshes and no significant difference in results or convergence 
behavior was observed.

\ms \noi
7.1.1 \ub{Constant $k$}
\par
\begin{figure}[h!]
	\centering
	\includegraphics[width=0.8\textwidth]{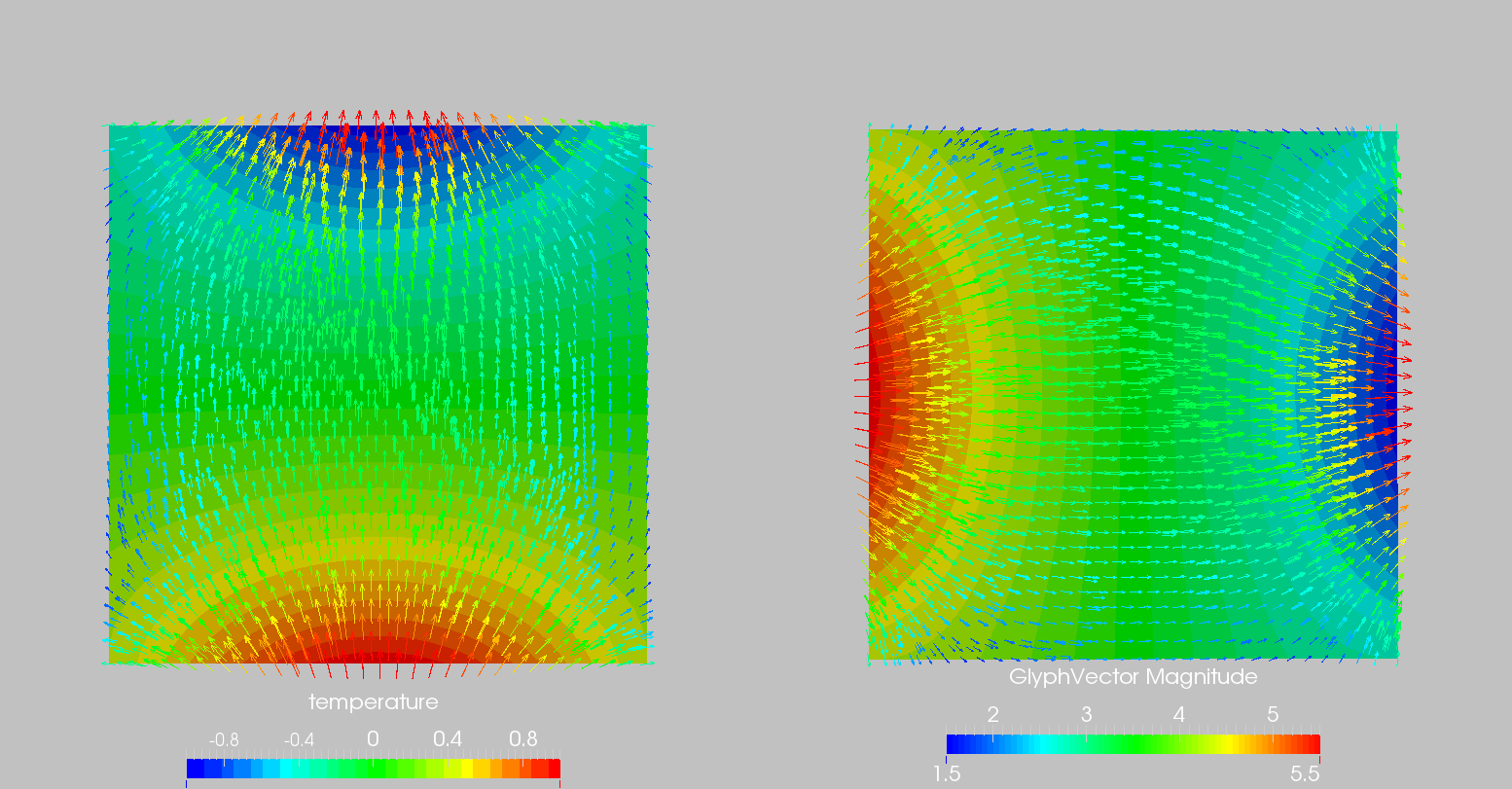}
	\caption{\label{f:square1}Square (a), (b) - constant $k$: Temperatures and Fluxes for Measurements 1,2.}
\end{figure}

\begin{figure}[h!]
	\centering
	\includegraphics[width=0.45\textwidth]{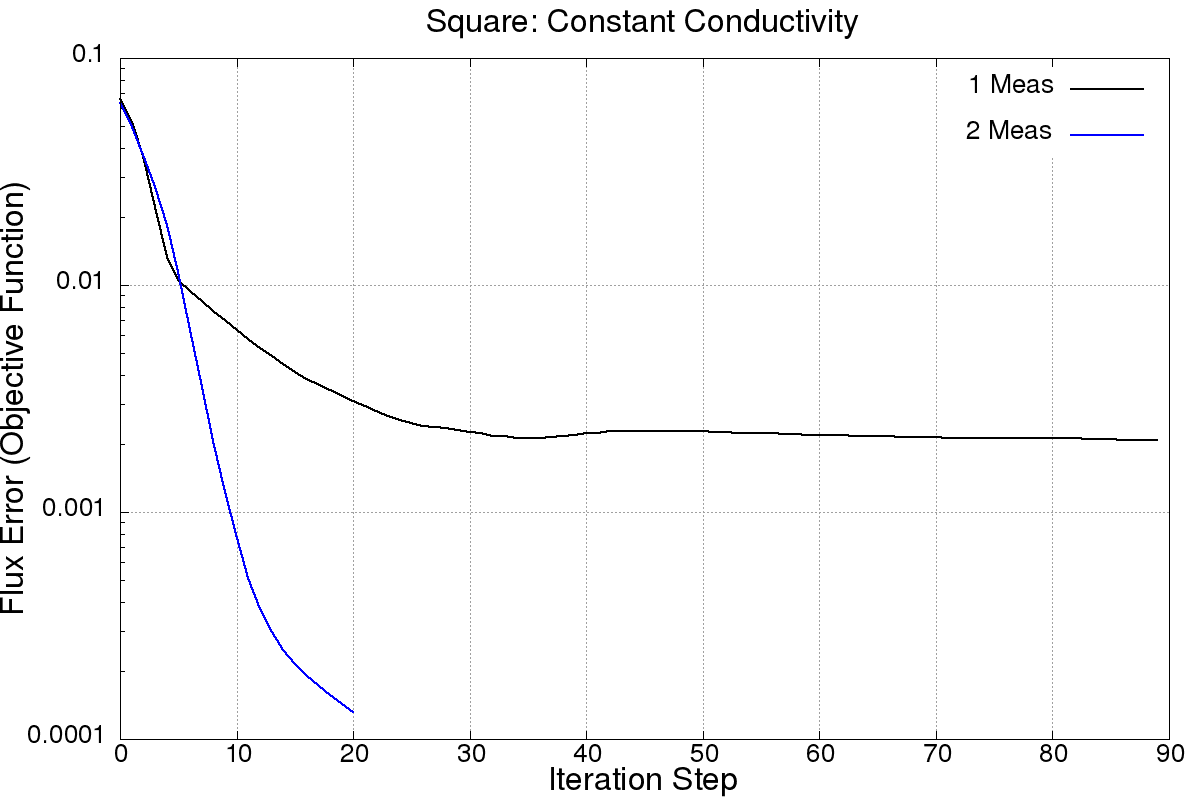}
	\includegraphics[width=0.45\textwidth]{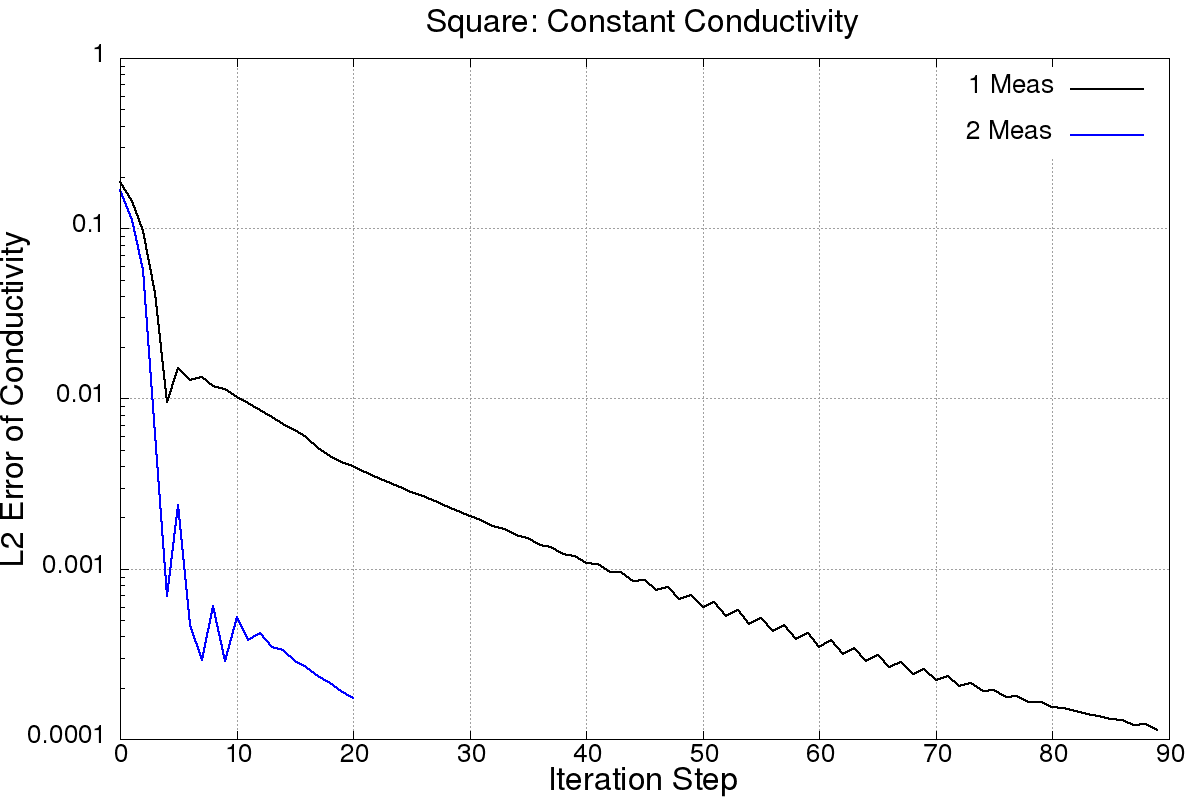}
	\caption{\label{f:square2}Square - constant $k$: Errors in Fluxes (left) and Conductivity (right).}
\end{figure}

 \noi
For this first test the desired (or target) conductivity $k$ 
is set to $k(\xvec)=2$. The temperature is set in a way that
is similar to the application of electric potentials in
Electrical Impedance Tomography \cite{DIsaacson_1986a,
lee1988determining,ZQChen_1990a,RVKohn_AMcKenney_1990a,
MCheney_DIsaacson_JCNewell_1999a,RHarikumar_RPrabu_SRaghavan_2013a,
BHarrach_JRieger_2019a,VChitturi_NFarrukh_2019a} via a set 
of `hot' and `cold' `sources' and `sinks' of the form:
\[
	u(\xvec) = \sum_i S_i(\xvec, \xvec_i, r_i, a_i)  \, ,
\]	
with
\[
	S_i(\xvec, \xvec_i, r_i, a_i) = 
		a_i e^{ {{(\xvec - \xvec_i)^2} \over {r_i^2}} }  \, .
\]
One source and one sink is used for each `measurement' as follows:
\begin{enumerate}
\item[-] First `Measurement' Parameters:
	\begin{enumerate}
		\item[-] $S_1: \xvec_1=(0.5,0.0)~,~r_1=0.5~,~a_1=1.0$ 
		\item[-] $S_2: \xvec_2=(0.5,1.0)~,~r_2=0.5~,~a_2=-1.0$ 
	\end{enumerate}		
\item[-] Second `Measurement' Parameters:
	\begin{enumerate}
		\item[-] $S_1: \xvec_1=(0.0,0.5)~,~r_1=0.5~,~a_1=1.0$ 
		\item[-] $S_2: \xvec_2=(1.0,0.5)~,~r_2=0.5~,~a_2=-1.0$ 
	\end{enumerate}		
\end{enumerate}
The case is first run with one and then with two measurements. 
Figures~\ref{f:square1} (a), (b) 
show the target temperature field and fluxes
field obtained after 10 steps (at this point the conductivity
is very close to the target conductivity and the difference in
temperature and fluxes between step 10 and subsequent steps negligible). 
The errors measured for the fluxes, as well as the L2 error of the 
conductivity are shown in Figures~\ref{f:square2} (left) and 
(right), respectively.  Noticeable in this result is the much 
faster rate of convergence when two measurements are used. 

\ms \noi
7.1.2 \ub{Linear $k$}
\par
\begin{figure}[h!]
	\centering
	\includegraphics[width=0.8\textwidth]{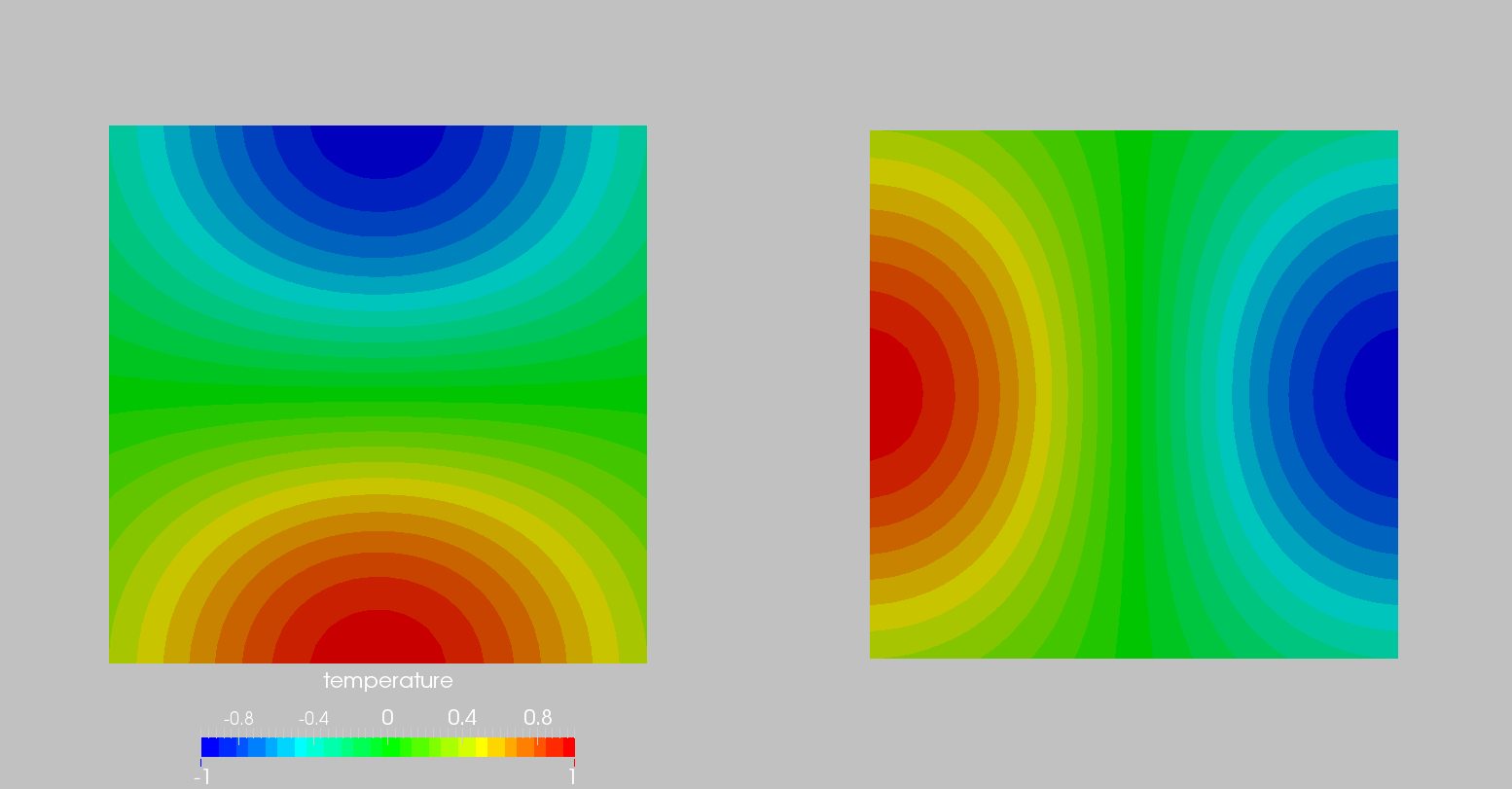}
	\includegraphics[width=0.8\textwidth]{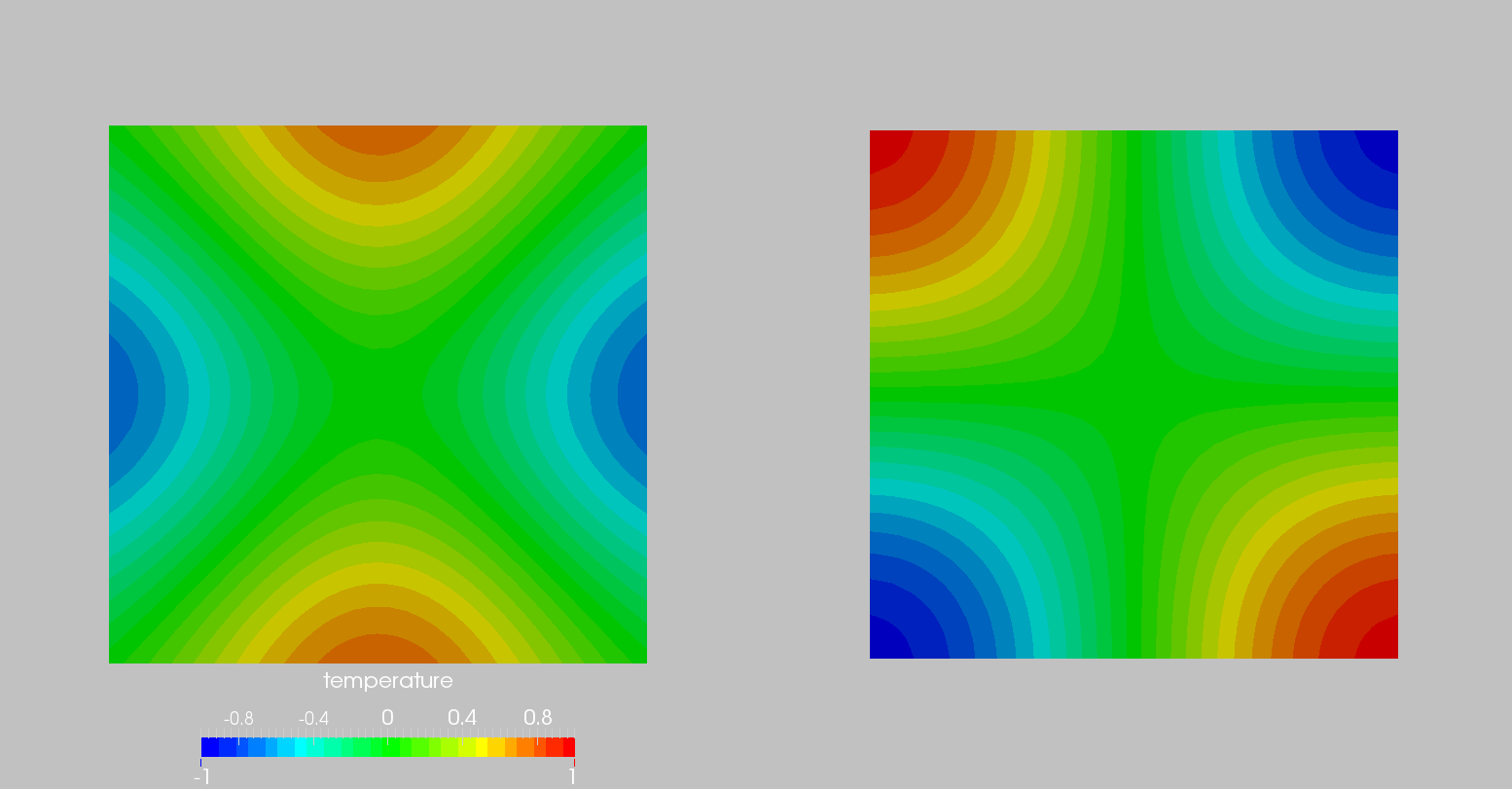}
	\caption{\label{f:squareLin}Square - linear $k$: Temperatures for Measurements 1-4.}
\end{figure}

\begin{figure}[h!]
	\centering
	\includegraphics[width=0.8\textwidth]{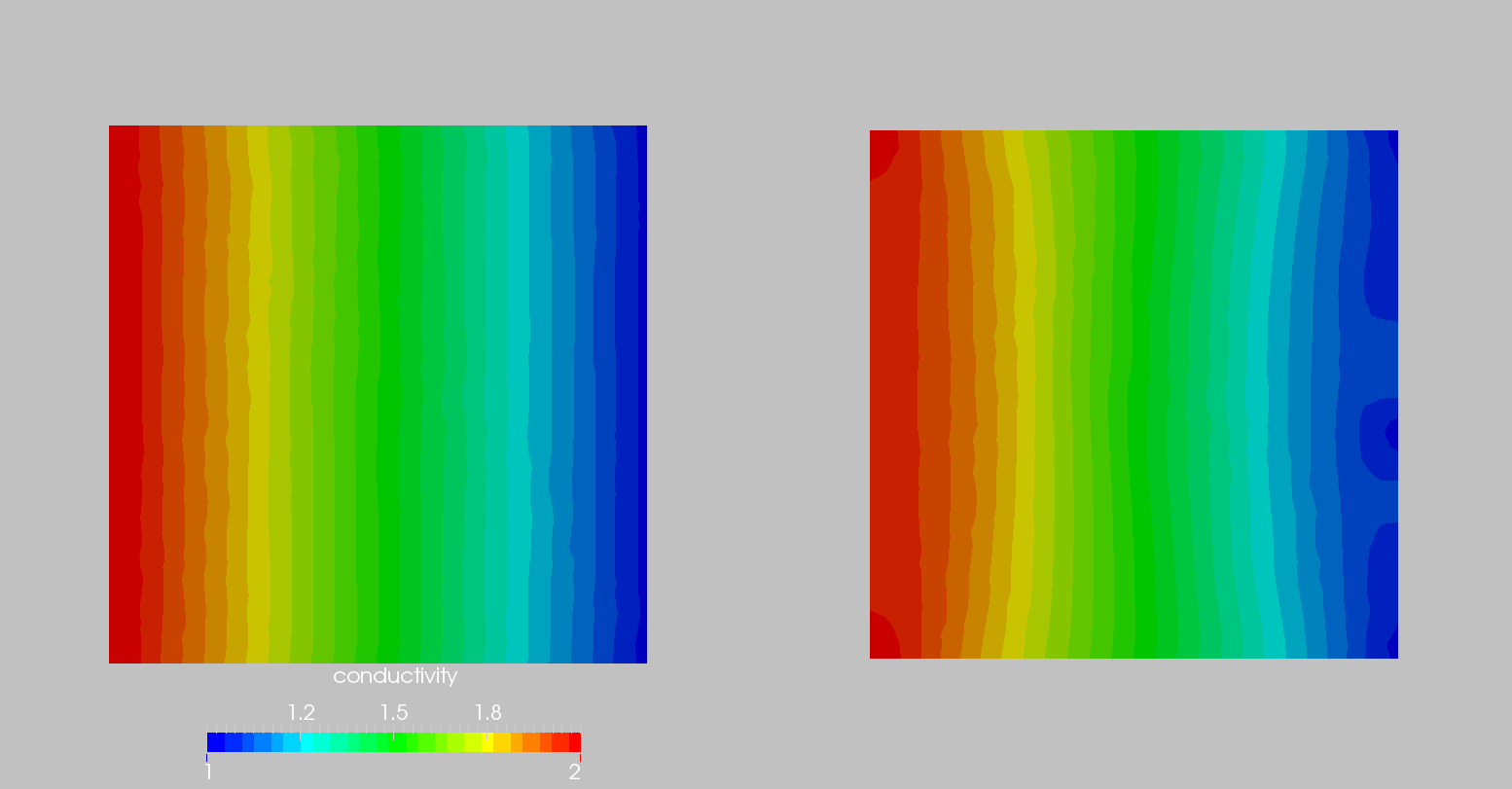}
	\caption{\label{f:squareLin2}Square - linear $k$: Target (Left) and Recovered (Right) Conductivity
with 4 Measurements.}
\end{figure}

\begin{figure}[h!]
	\centering
	\includegraphics[width=0.45\textwidth]{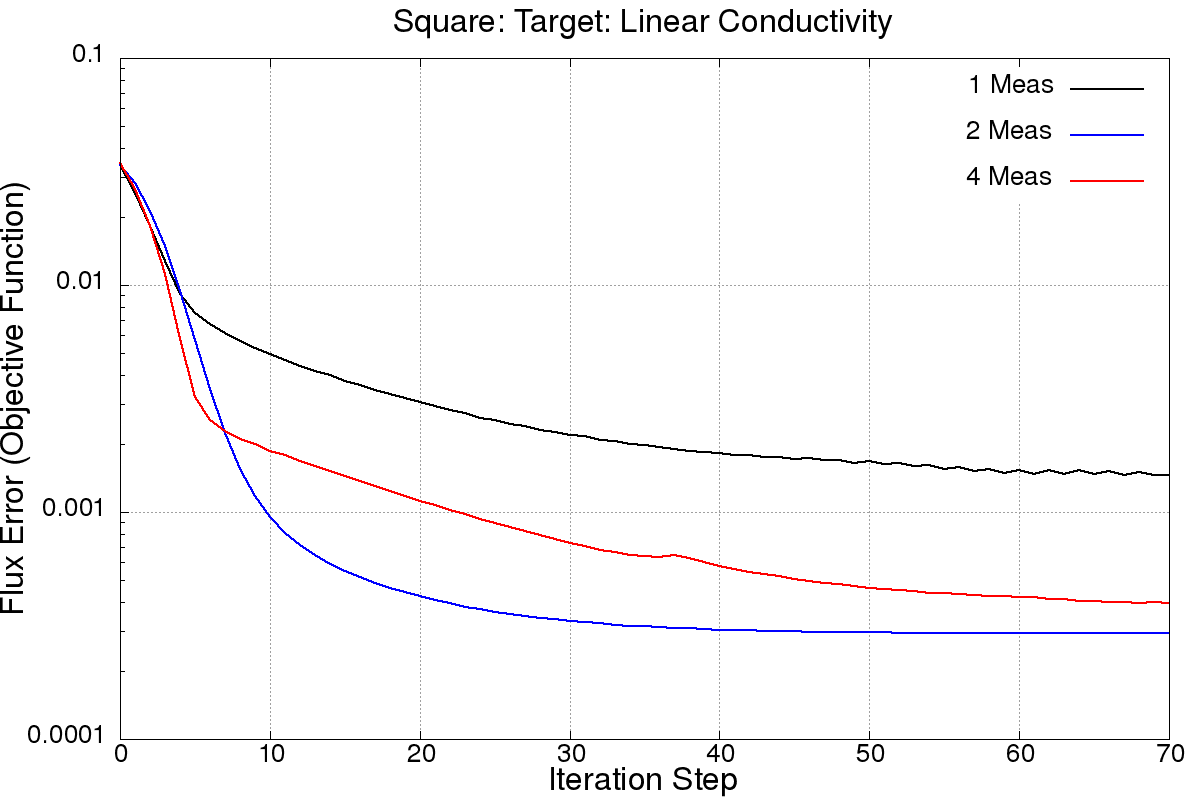}
	\includegraphics[width=0.45\textwidth]{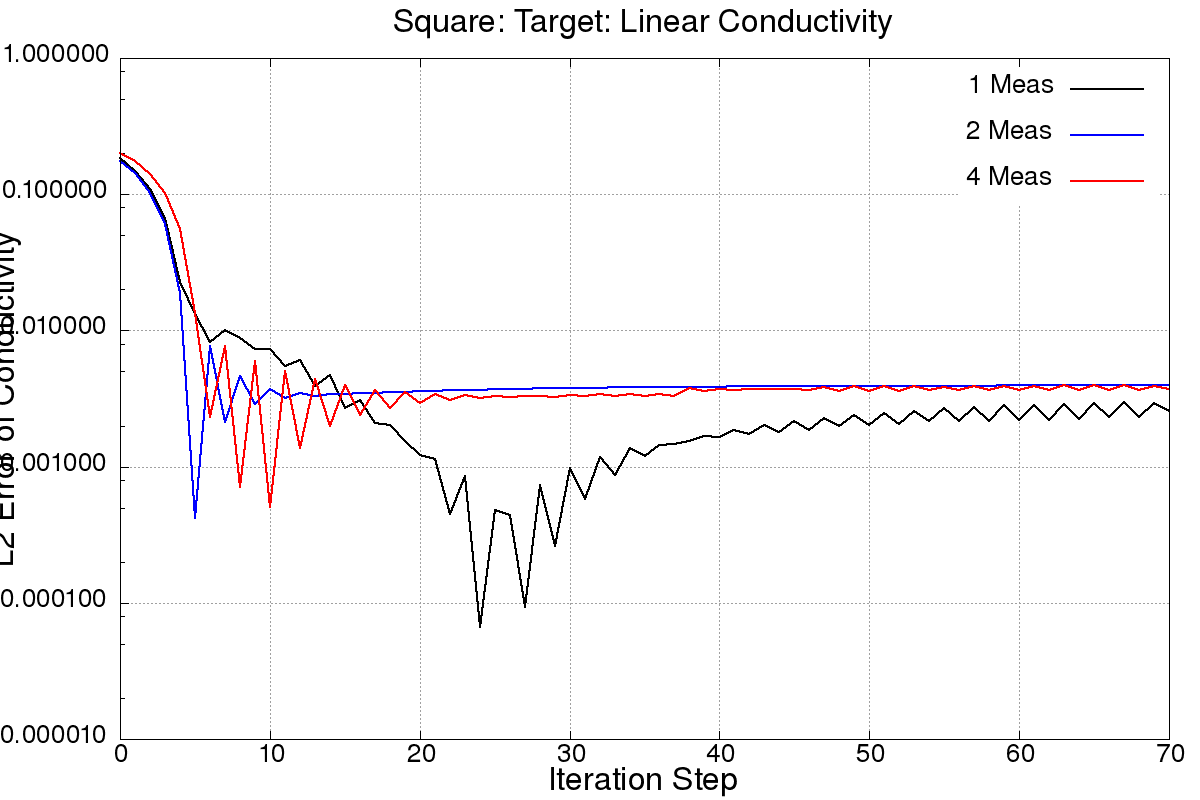}
	\caption{\label{f:squareLin3}Square - linear $k$: Errors in Fluxes (left) and Conductivity (right).}
\end{figure}

 \noi
The next test considers as a target a linear conductivity of the form:
\[
	k(\xvec) = 2 - x \, .
\]
The up to four temperature fields prescribed (i.e. measurements) are
shown in Figure~\ref{f:squareLin}. 
The conductivity $k$ in each element was set
constant, based on the element centroid (average of the 4 nodes).
Figures~\ref{f:squareLin2} show the target conductivity (left) 
and the conductivity
obtained with four measurements (right). Note that because the graphics 
first average to point values and then perform the plane cut for display
even the target (exact) conductivity appears noisy. Nevertheless, a
difference between the target and recovered conductivity may be observed.
The errors measured for the fluxes, as well as the L2 error of the 
conductivity are shown in Figures~\ref{f:squareLin3} 
(left) and (right), respectively. As before, a single measurement 
seems to lead to slower convergence. It also appears that the number of
measurements at some point `saturates': 4 measurements do not yield a
better or faster result than 2 measurements.

\ms \noi
7.1.3 \ub{Gaussian $k$}
\par \noi
In this case the target distribution for $k$ is given by: 
\[
	k(\xvec)=1+4 e^{\big({{|\xvec-\xvec_0|} \over r_0}\big)^2} \, , \quad 
		\xvec_0=(0.5,0.5,0.025) \, , \quad r_0=0.2 \, . 
\]		

\begin{figure}[h!]
	\centering
	\includegraphics[width=0.8\textwidth]{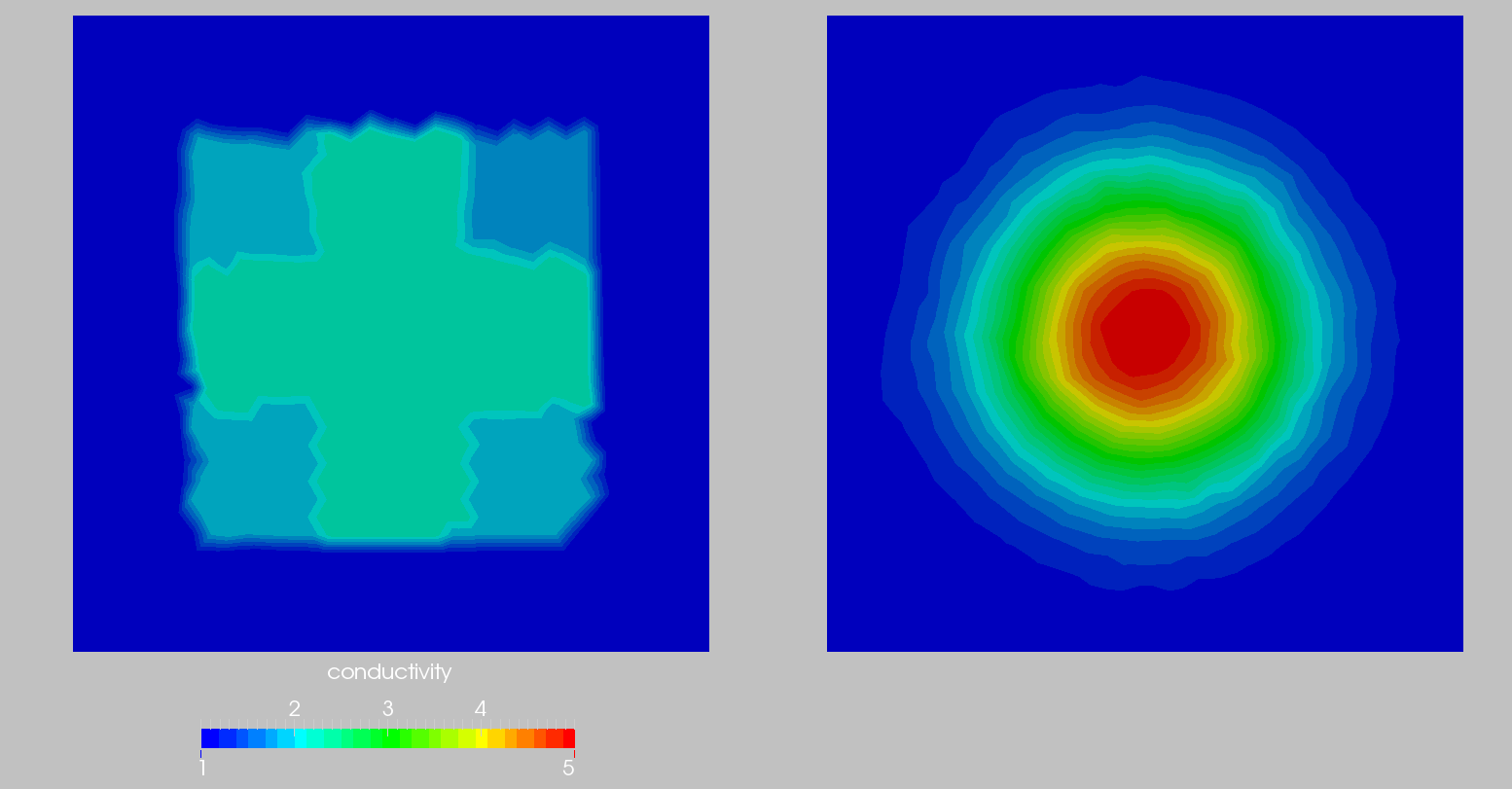}
	\caption{\label{f:squareGauss}Square - Gaussian $k$: Conductivity Distribution: Left: 25~DOFs,
Right: Target.}
\end{figure}

\begin{figure}[h!]
	\centering
	\includegraphics[width=0.8\textwidth]{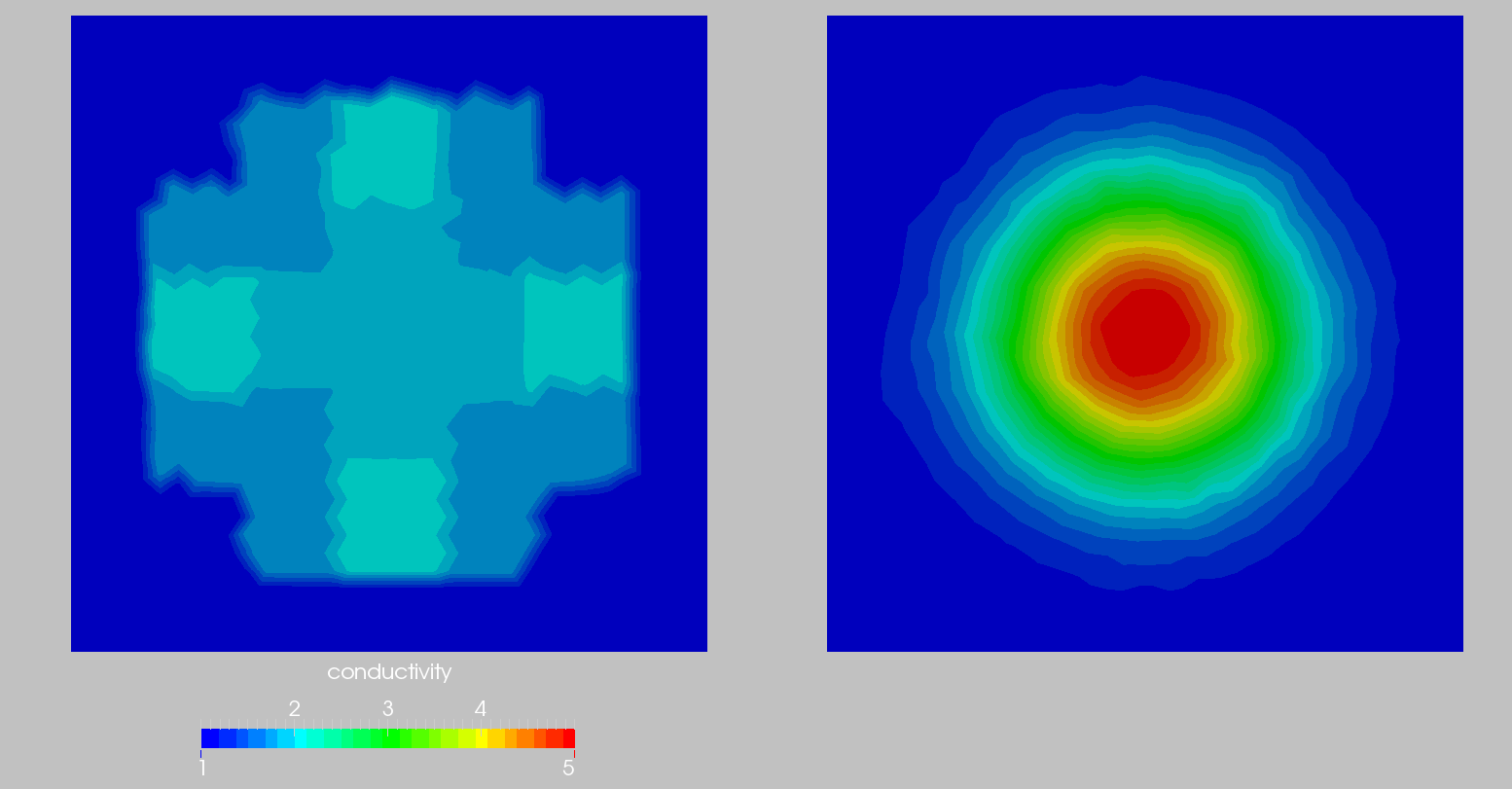}
	\caption{\label{f:squareGauss2}Square - Gaussian $k$: Conductivity Distribution: Left: 49~DOFs,
Right: Target.}
\end{figure}

\begin{figure}[h!]
	\centering
	\includegraphics[width=0.8\textwidth]{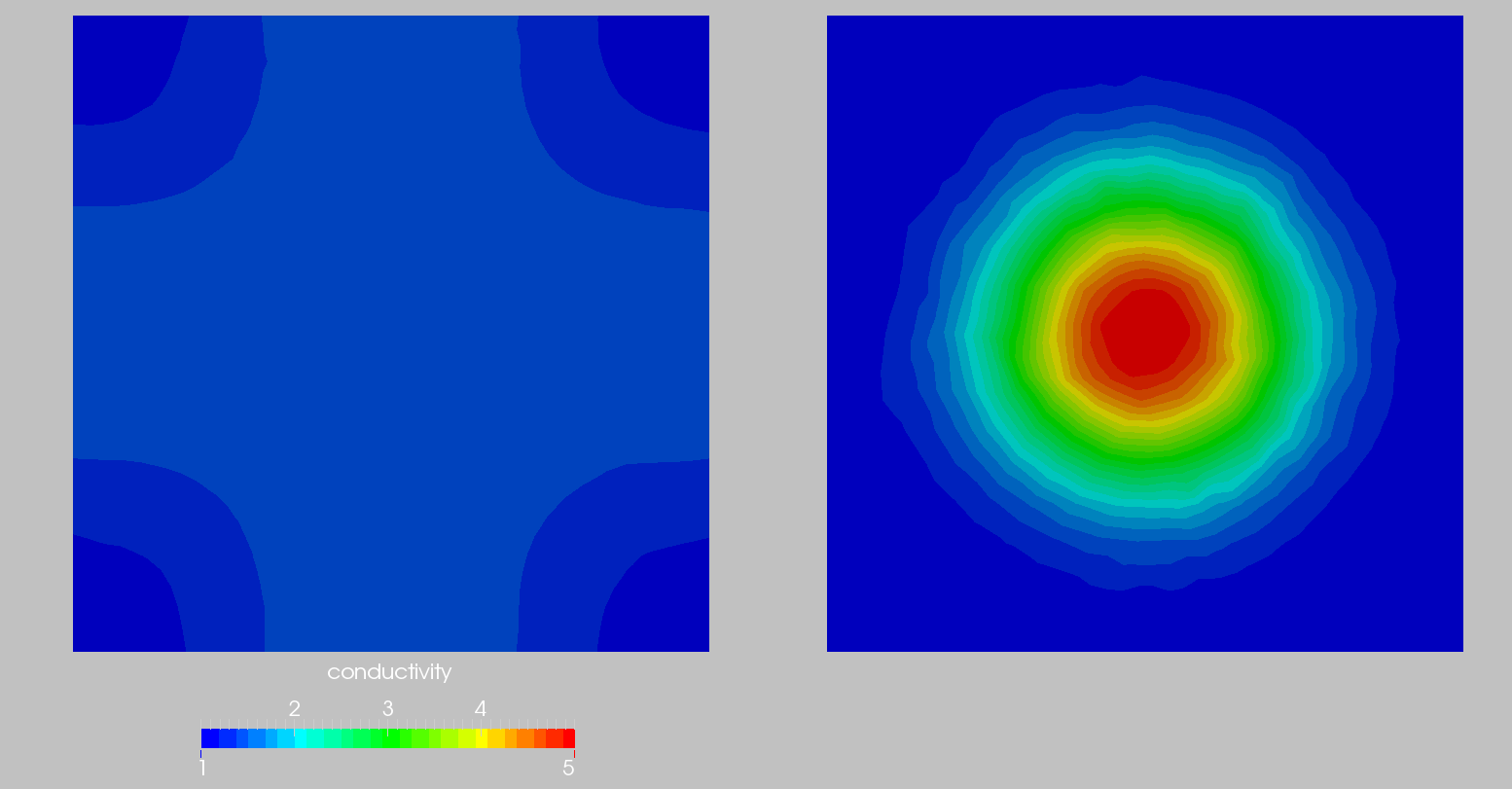}
	\caption{\label{f:squareGauss3}Square - Gaussian $k$: Conductivity Distribution: 
Left: Adjoint-Based, Right: Target.}
\end{figure}

\begin{figure}[h!]
	\centering
	\includegraphics[width=0.45\textwidth]{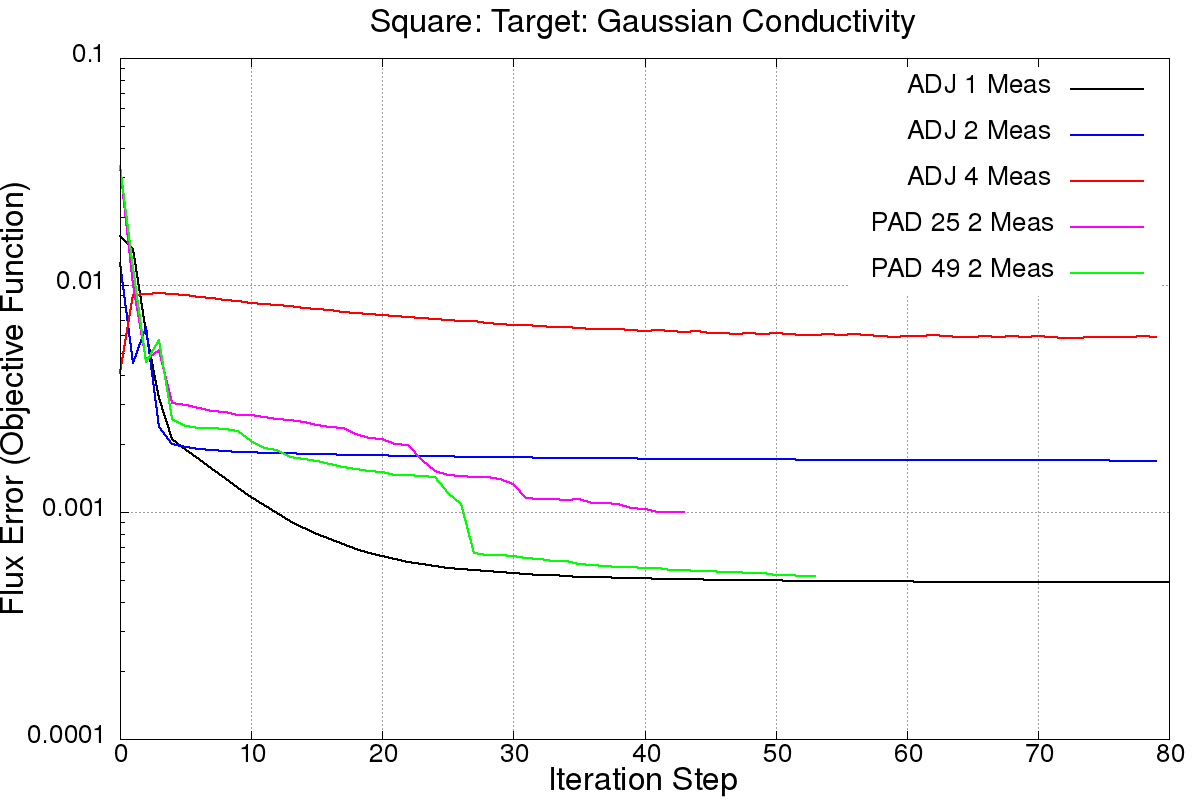}
	\includegraphics[width=0.45\textwidth]{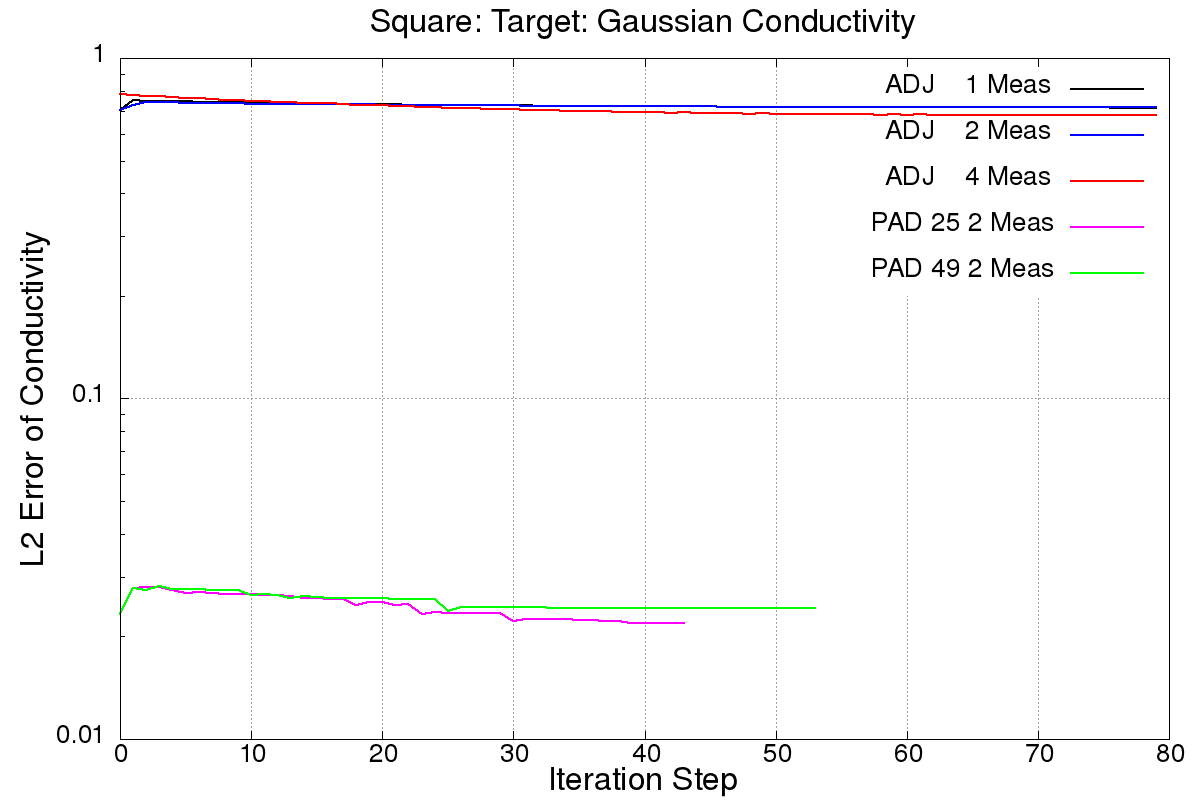}
	\caption{\label{f:squareGauss4}Square - Gaussian $k$: Errors in Fluxes (left) and Conductivity (right).}
\end{figure}

\noi
The same 4 `measurement functions' as before are used for this case.
For the optimization via finite differences, the distribution
of $k$ was given by constant regions spaced in a cartesian lattice in
$x,y$. Figures~\ref{f:squareGauss}, \ref{f:squareGauss2}, 
\ref{f:squareGauss3} show the results obtained using 25 and 49
degrees of freedom for the finite difference gradient and
projected adjoint gradient optimization, as well as the adjoint. 
The decrease in objective function for these
cases is shown in Figure~\ref{f:squareGauss4} (left), and the L2 
error of the conductivity in Figure~\ref{f:squareGauss4} (right). 
Note that the optimization procedure as
such works, recovering well the specified normal flux distribution on
the boundary. Note the slow convergence of the adjoint, which may
be attributable to the much larger degrees of freedom. Note also
that the spatial distribution for $k$, while resembling the target 
distribution, still has considerable differences. Seemingly, the
lower the degrees of freedom, the better the result (!).

\ms \noi
7.1.4 \ub{Disk-Like $k$}
\par
\begin{figure}[h!]
	\centering
	\includegraphics[width=0.8\textwidth]{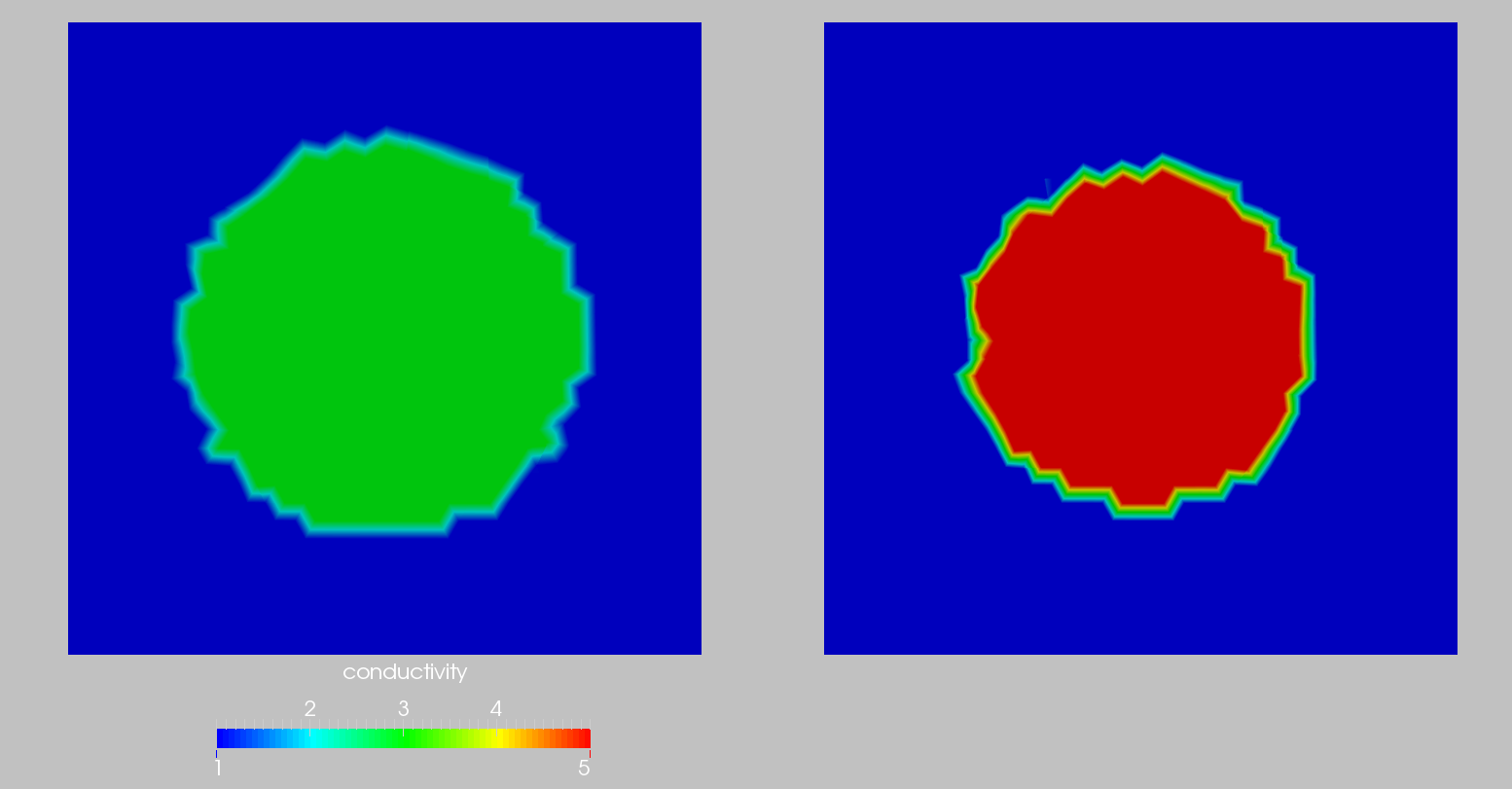}
	\caption{\label{f:squarePuck} Square - Disk-Like $k$: Conductivity Distribution:
Left: Converged Result, Right: Target}
\end{figure}

\begin{figure}[h!]
	\centering
	\includegraphics[width=0.8\textwidth]{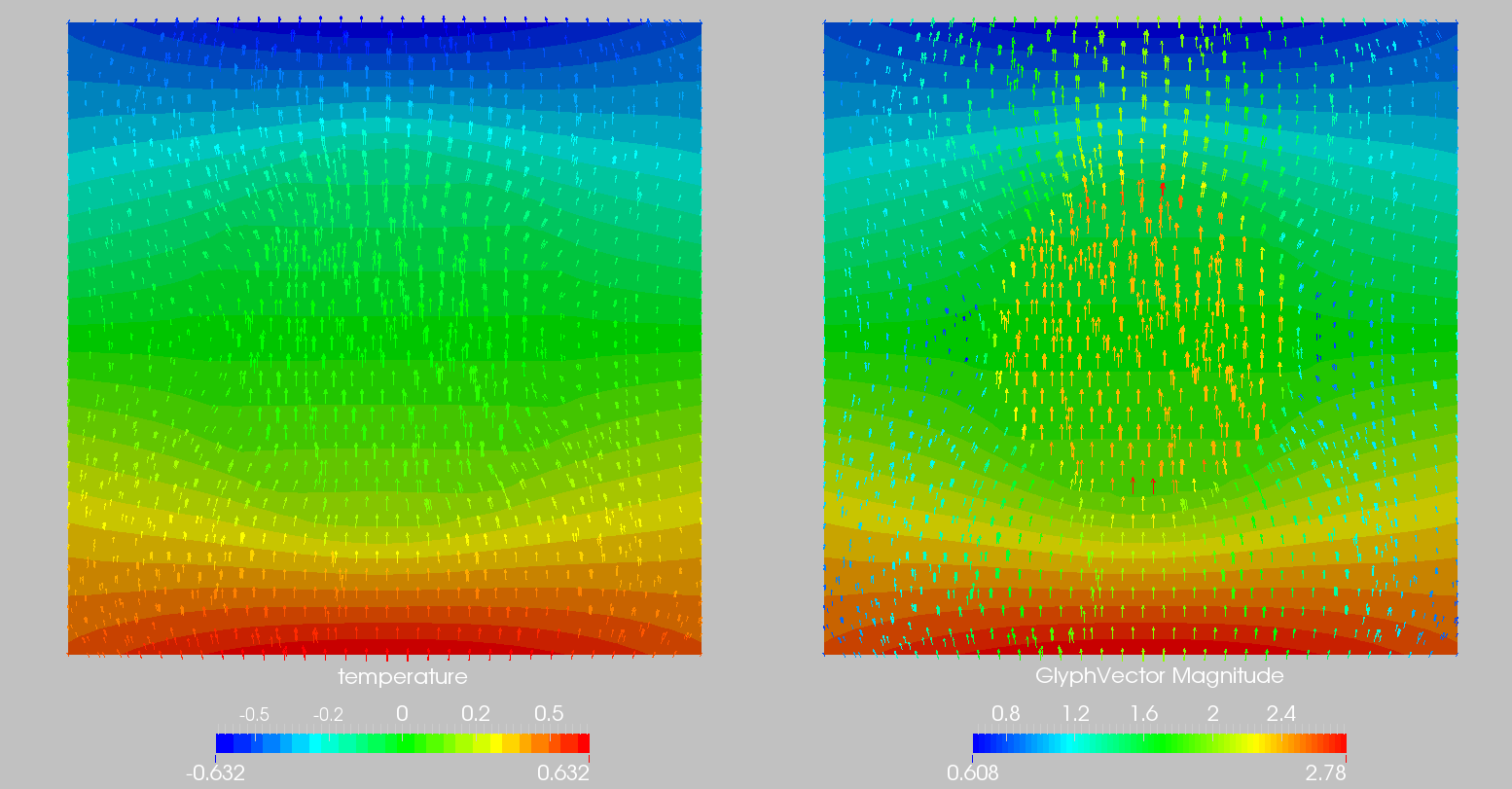}
	\caption{\label{f:squarePuck2} Square - Disk-Like $k$: Temperature and Flux Distribution:
Left: Converged Result, Right: Target}
\end{figure}

\begin{figure}[h!]
	\centering
	\includegraphics[width=0.5\textwidth]{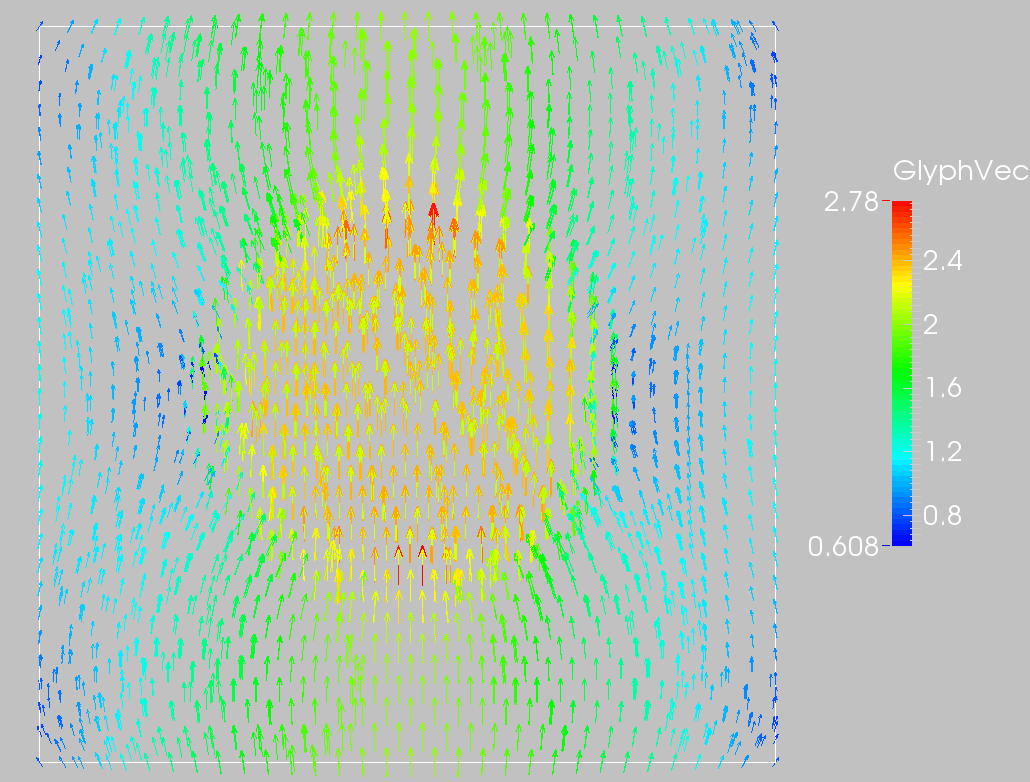}
	\caption{\label{f:squarePuck3} Square - Disk-Like $k$: Comparison of Flux Distribution for 
Converged Result and Target}
\end{figure}

\begin{figure}[h!]
	\centering
	\includegraphics[width=0.45\textwidth]{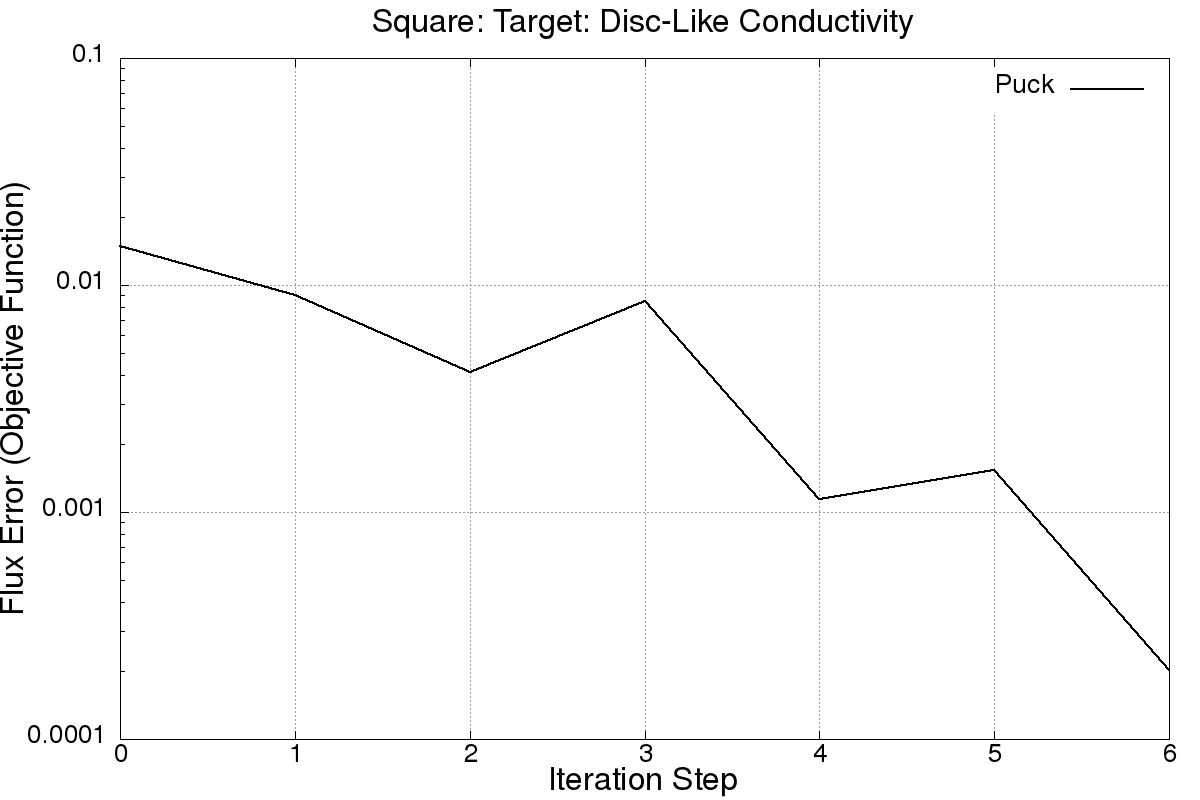}
	\includegraphics[width=0.45\textwidth]{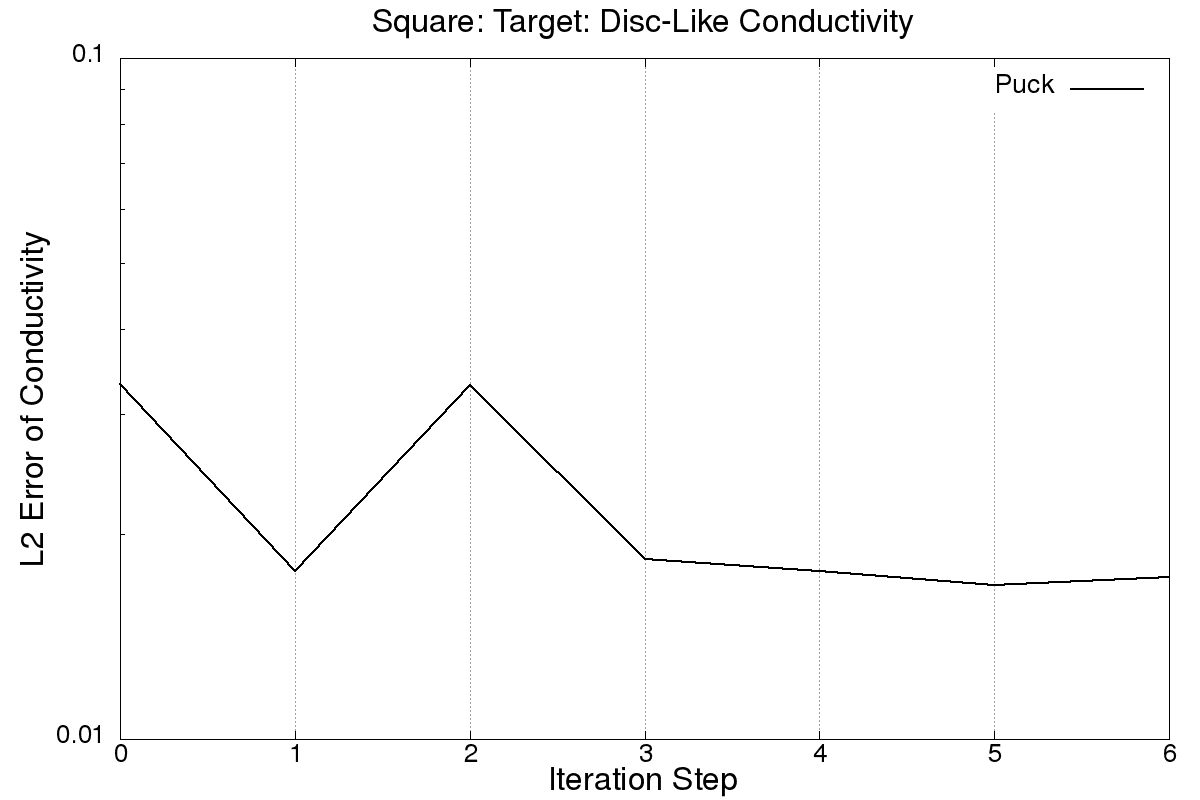}	
	\caption{\label{f:squarePuck4} Square - Disk-Like $k$: Errors in Fluxes (left) and Conductivity (right)}
\end{figure}

 \noi
Intrigued by the previous result, a case with very few design
parameters has been chosen. The (discontinuous) target distribution 
for $k$ of the form:
\[ 
\begin{aligned}
k(\xvec)=k_{disk}=5 ~~ &\forall \ |\xvec-\xvec_0| \le r_0 \, , \quad 
   k(\xvec)=k_{exte}=1 ~~ \forall \ |\xvec-\xvec_0| > r_0 ~~,~~ \\
   &\xvec_0=(0.5,0.5,0.025) ~~,~~ r_0=0.25 ~~. 
\end{aligned}   
\]   
Two `measurement functions' were used for this case.
The case was run with only 4 design parameters:
$x_0, y_0, r_0$ and $k_{disc}$. The initial values were set to:
$x_0=y_0=0.25, r_0=0.1, k_{disc}=2$, i.e. relatively far away from
the target parameters. The results obtained are shown in 
Figures~\ref{f:squarePuck},
\ref{f:squarePuck2}, \ref{f:squarePuck3}. 
The decrease in objective function for these
cases is shown in Figure~\ref{f:squarePuck4} (left), and the L2 error 
of the conductivity in Figure~\ref{f:squarePuck4} (right). Note the 
fast convergence and the correct position
for the target conductivity. The radius obtained is a bit larger, which
implies that the conductivity obtained $k_{disk}$ has to be lower.
However, the normal boundary fluxes are indistinguishable from one
another (see Figure~\ref{f:squarePuck3}).
We remark that the optimizer has converged to a relative design parameter 
range of $\epsilon_r=10^{-3}$, i.e. both $r_0$ and $k_{disk}$ will not 
change significantly even if $\epsilon_r$ is set lower.

\clearpage

\subsection{Cube}

This case is similar to the one shown before, but fully 3-D.
The domain dimensions are
$0 \le x \le 1, 0 \le y \le 1, 0 \le z \le 1$. 
The finite element mesh had more than 250,000 tetrahedral 
elements of uniform size.

\ms \noi
7.2.1 \ub{Gaussian $k$}
\par \noi
In this case the target distribution for $k$ is given by:
\[
	k(\xvec)=1+4 e^{\big({{|\xvec-\xvec_0|} \over r_0}\big)^2} \, \quad 
		\xvec_0=(0.5,0.5,0.05) \, \quad r_0=0.2 \, . 
\]		

\begin{figure}[h!]
	\centering
	\includegraphics[width=0.8\textwidth]{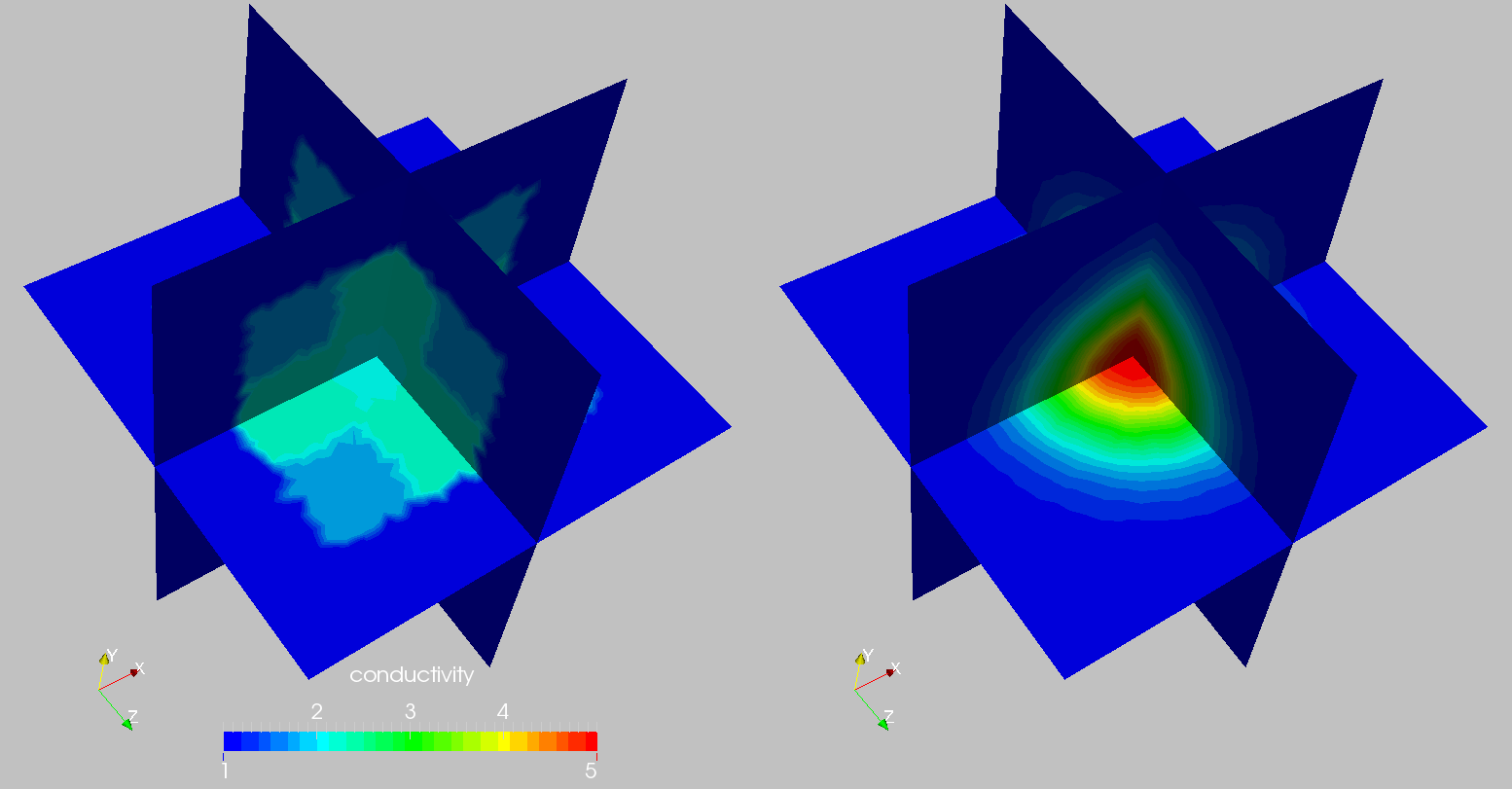}
	\caption{\label{f:cubeGauss} Cube - Gaussian $k$: Conductivity Distribution: Left: 125~DOFs,
Right: Target}
\end{figure}

\begin{figure}[h!]
	\centering
	\includegraphics[width=0.8\textwidth]{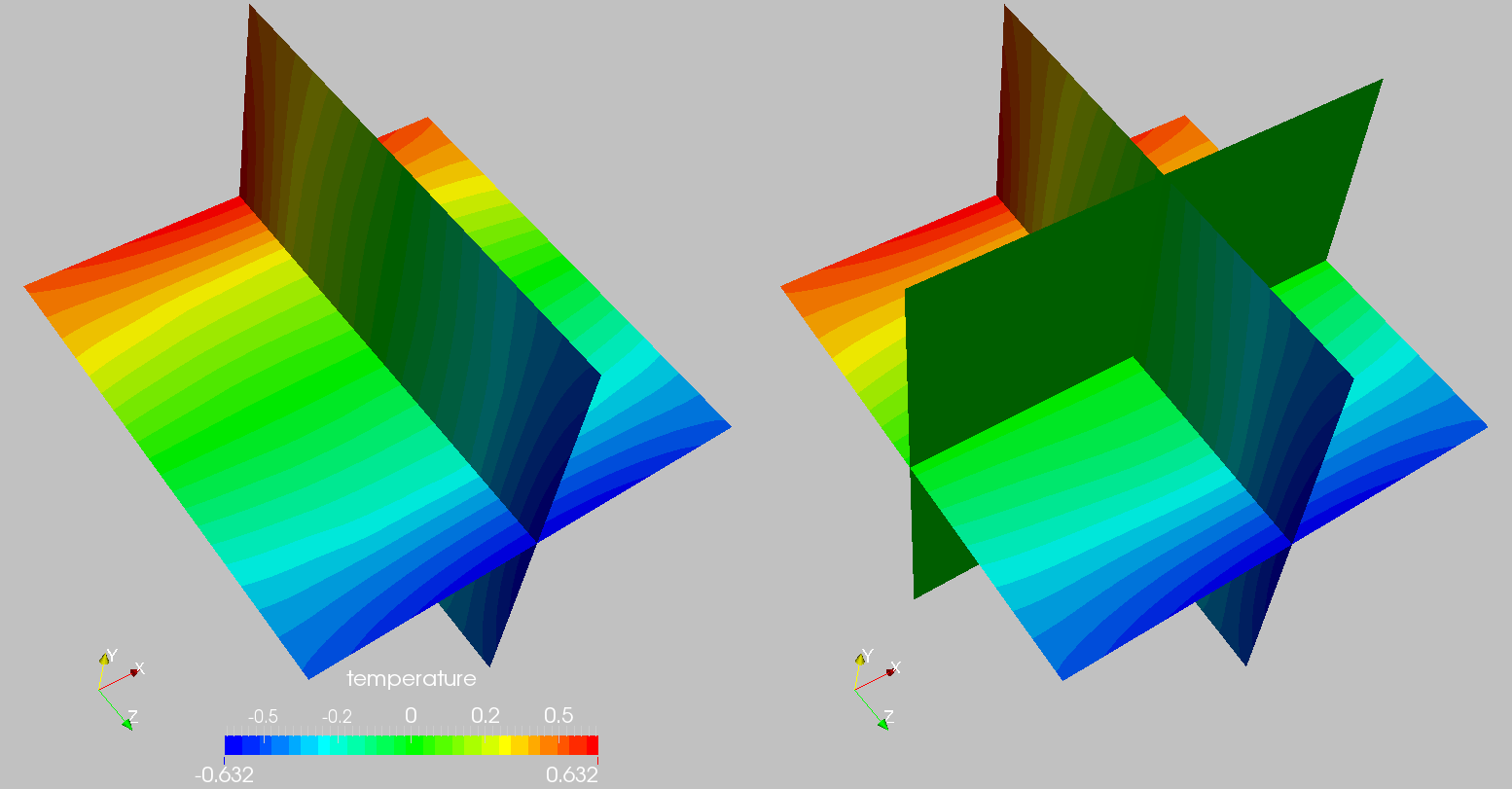}
	\caption{\label{f:cubeGauss2} Cube - Gaussian $k$: Temperature Distribution for Measurement 3: 
Left: 125~DOFs, Right: Target}
\end{figure}

\begin{figure}[h!]
	\centering
	\includegraphics[width=0.8\textwidth]{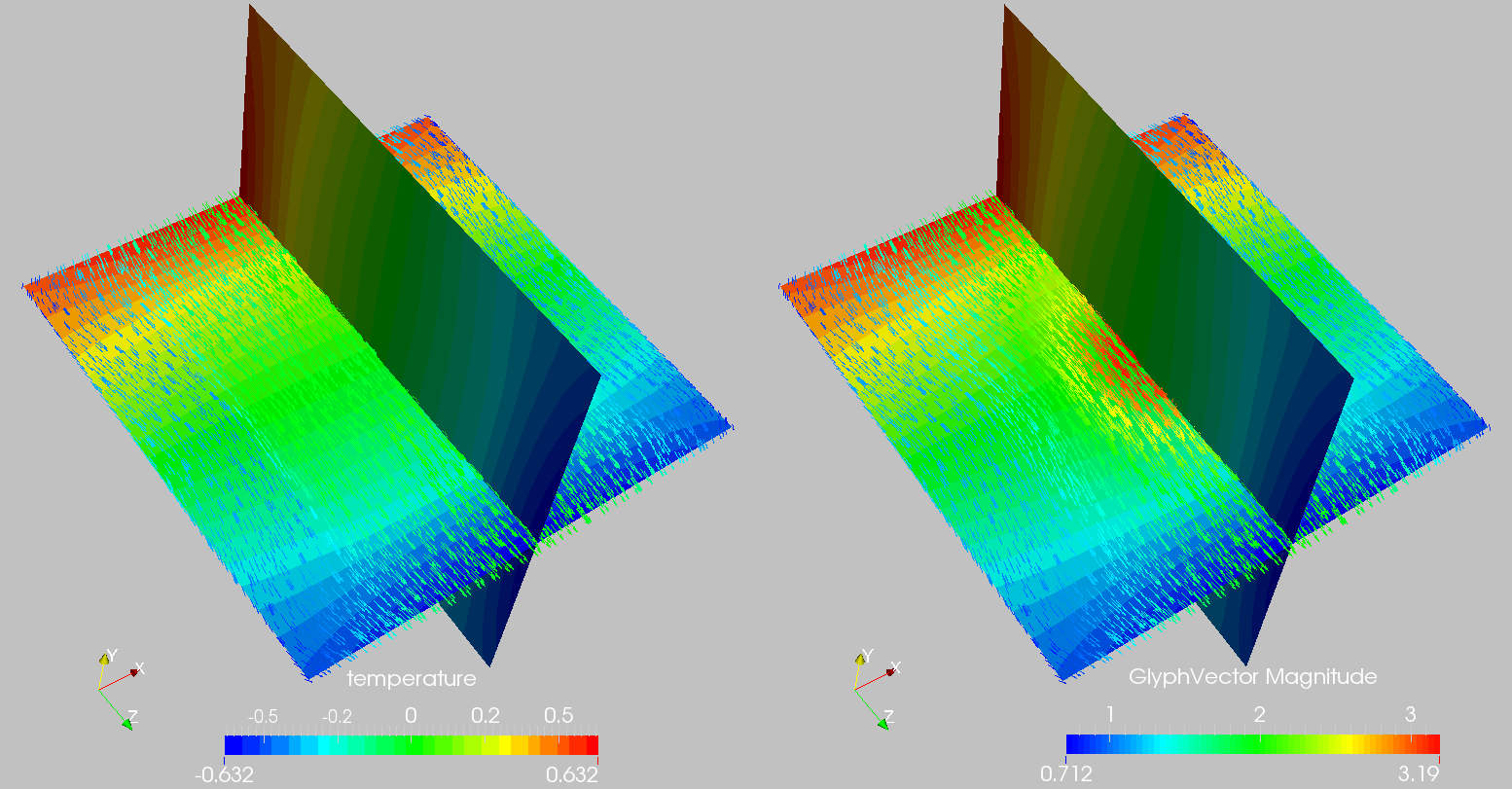}
	\caption{\label{f:cubeGauss3} Cube - Gaussian $k$: Flux Distribution for Measurement 3: 
Left: 125~DOFs, Right: Target}
\end{figure}

\begin{figure}[h!]
	\centering
	\includegraphics[width=0.8\textwidth]{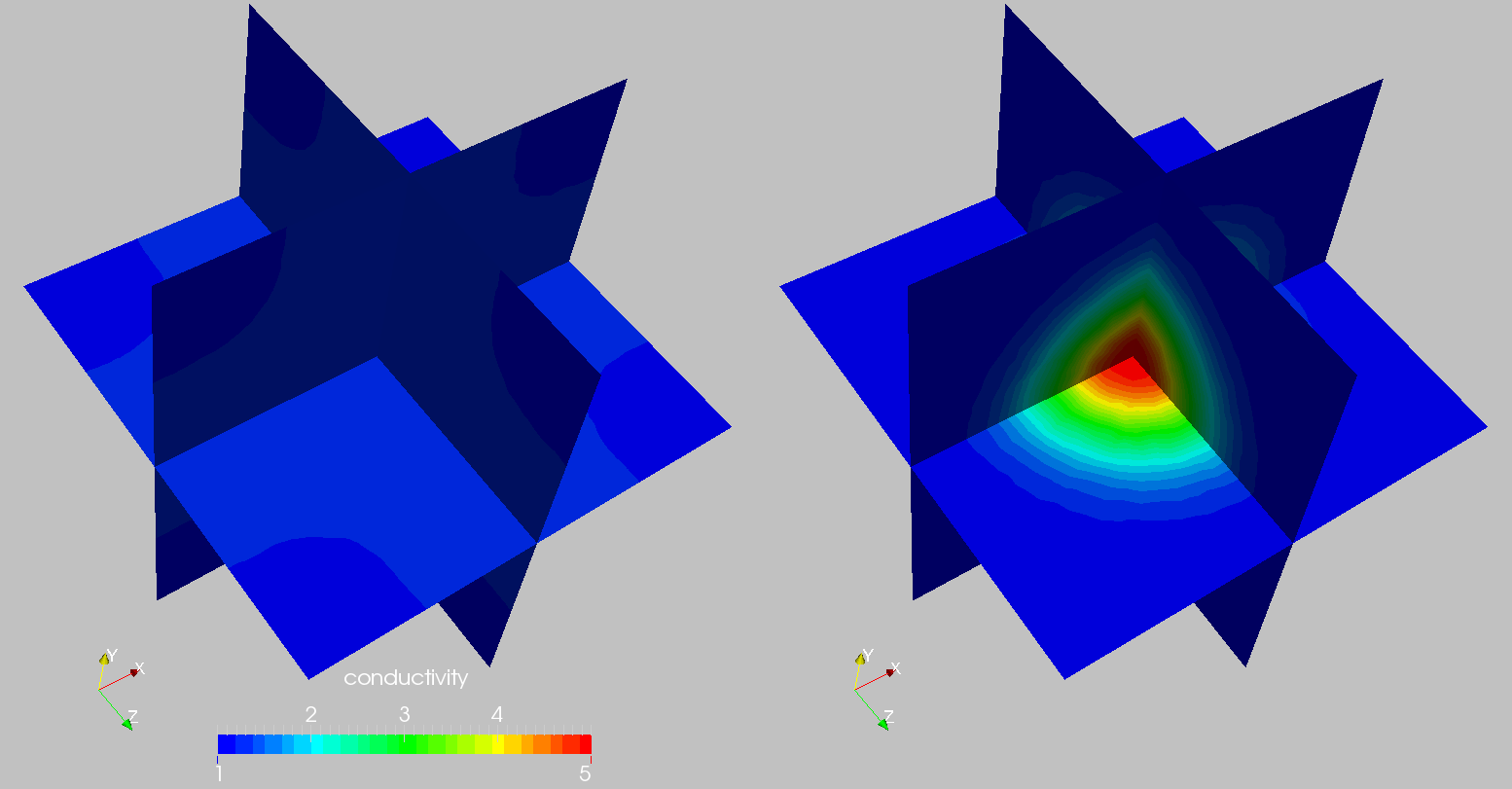}
	\caption{\label{f:cubeGauss4} Cube - Gaussian $k$: Conductivity Distribution: Left: 250K~DOFs,
Right: Target}
\end{figure}

\begin{figure}[h!]
        \centering
        \includegraphics[width=0.8\textwidth]{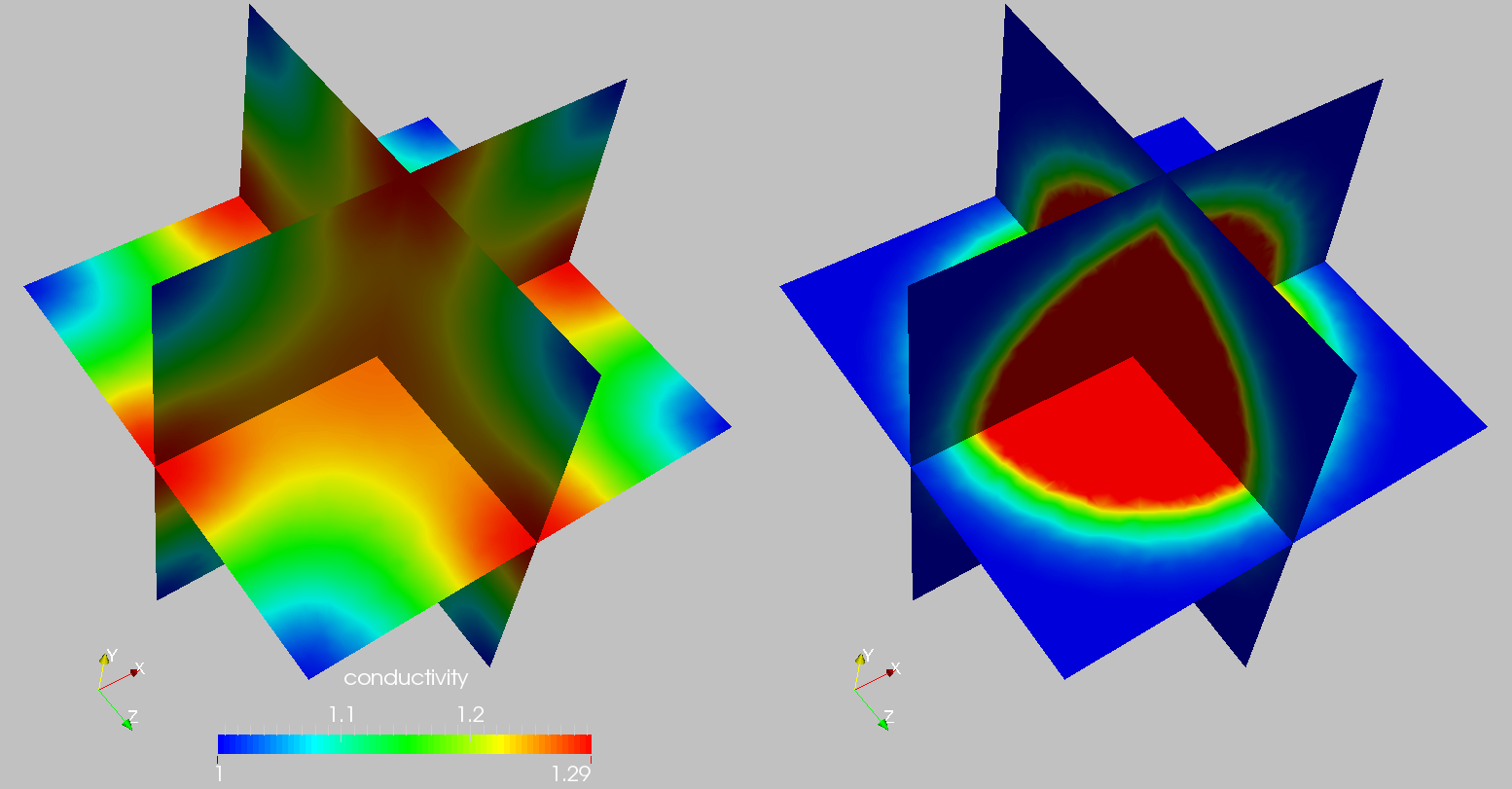}
        \caption{\label{f:cubeGauss4a} Cube - Gaussian $k$: Conductivity Distribution: Left: 250K~DOFs, Right: Target (Note: Re-Scaled for Visualization)}
\end{figure}

\begin{figure}[h!]
	\centering
	\includegraphics[width=0.8\textwidth]{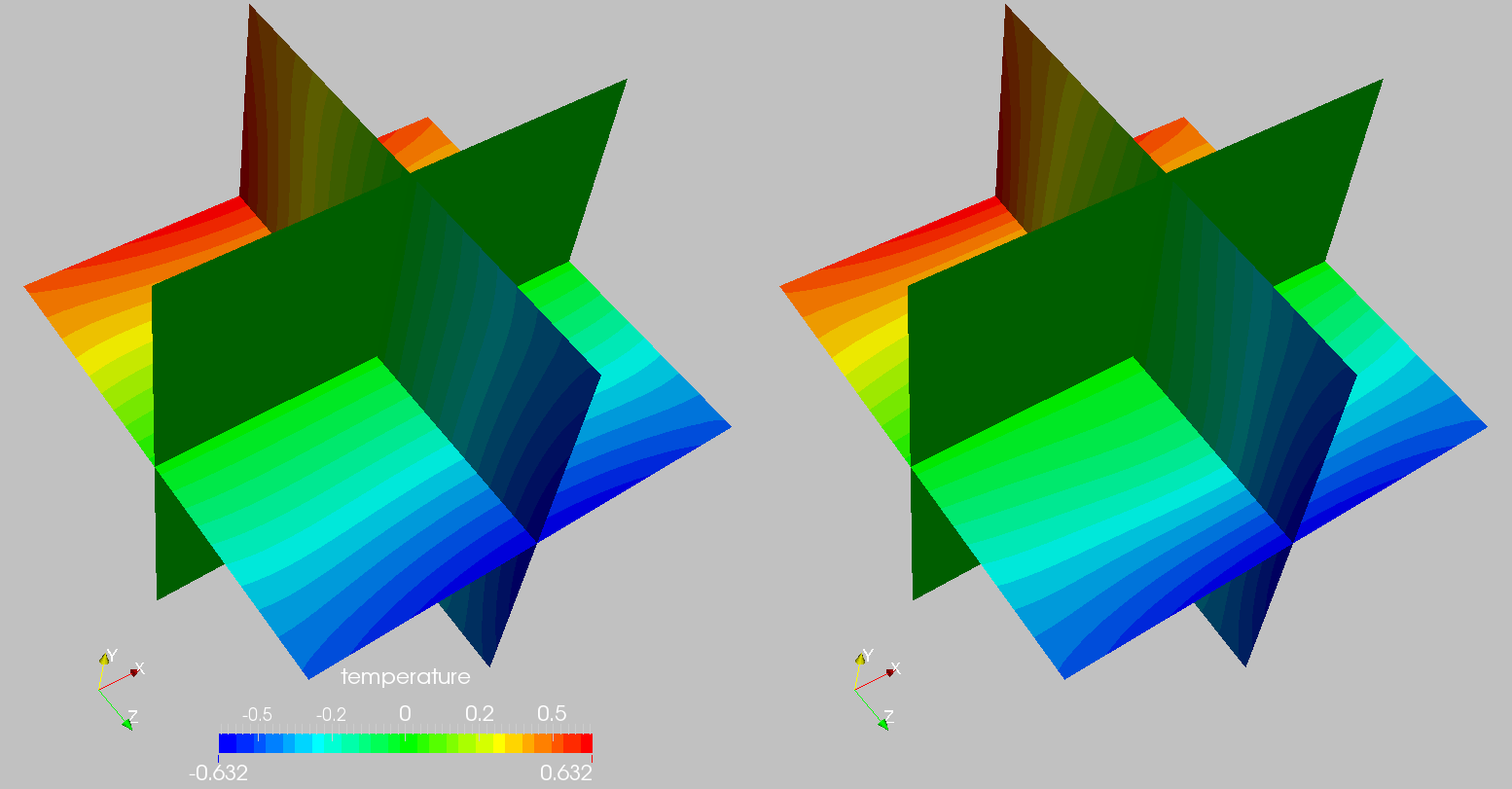}
	\caption{\label{f:cubeGauss5} Cube - Gaussian $k$: Temperature Distribution for Measurement 3:
Left: 250K~DOFs, Right: Target}
\end{figure}

\begin{figure}[h!]
	\centering
	\includegraphics[width=0.8\textwidth]{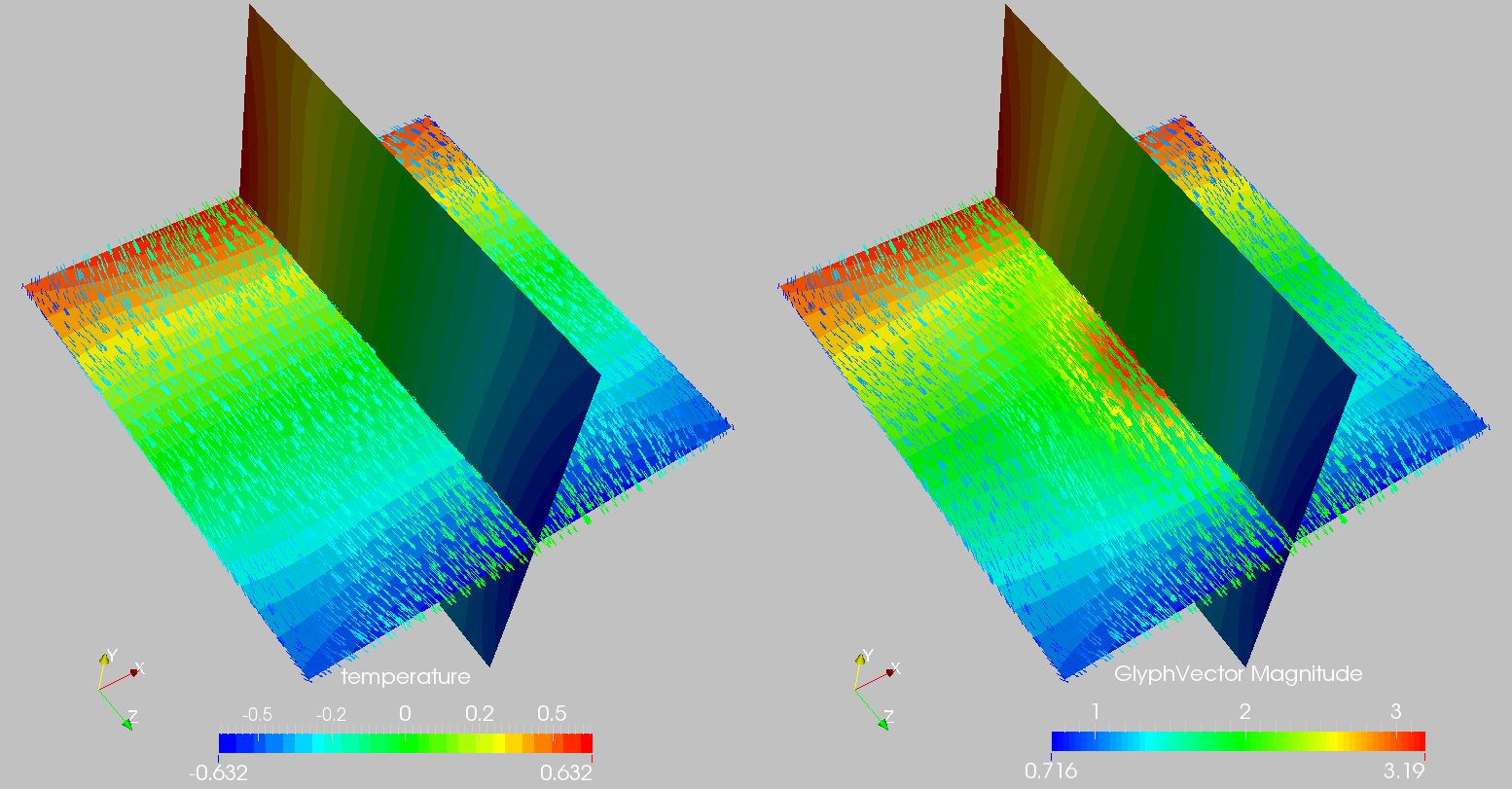}
	\caption{\label{f:cubeGauss6} Cube - Gaussian $k$: Flux Distribution for Measurement 3:
Left: 250K~DOFs, Right: Target}
\end{figure}

\begin{figure}[h!]
	\centering
	\includegraphics[width=0.5\textwidth]{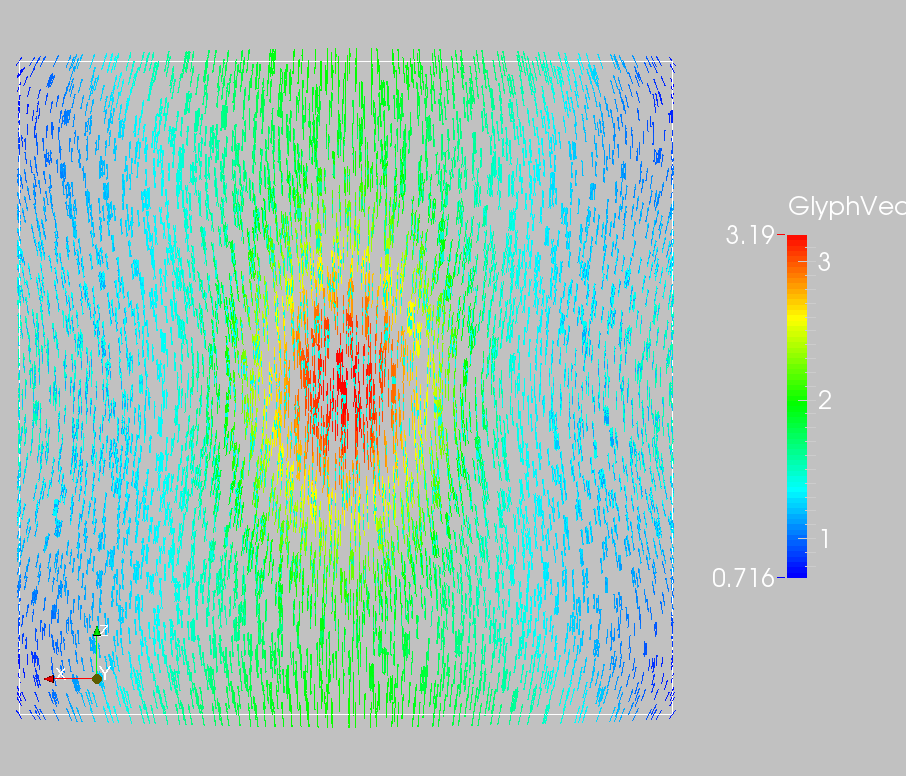}
	\caption{\label{f:cubeGauss7} Cube - Gaussian $k$: Comparison of Fluxes for Plane $y=0$}
\end{figure}

\begin{figure}[h!]
	\centering
	\includegraphics[width=0.45\textwidth]{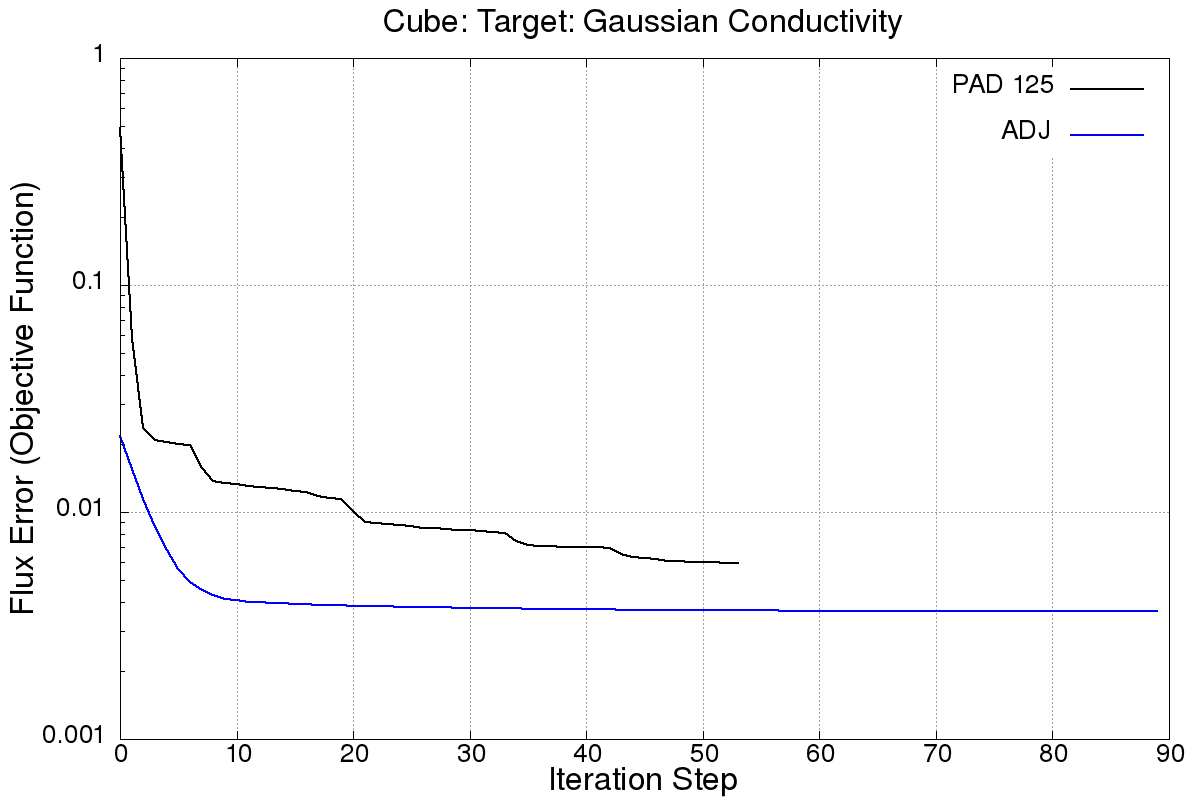}
	\includegraphics[width=0.45\textwidth]{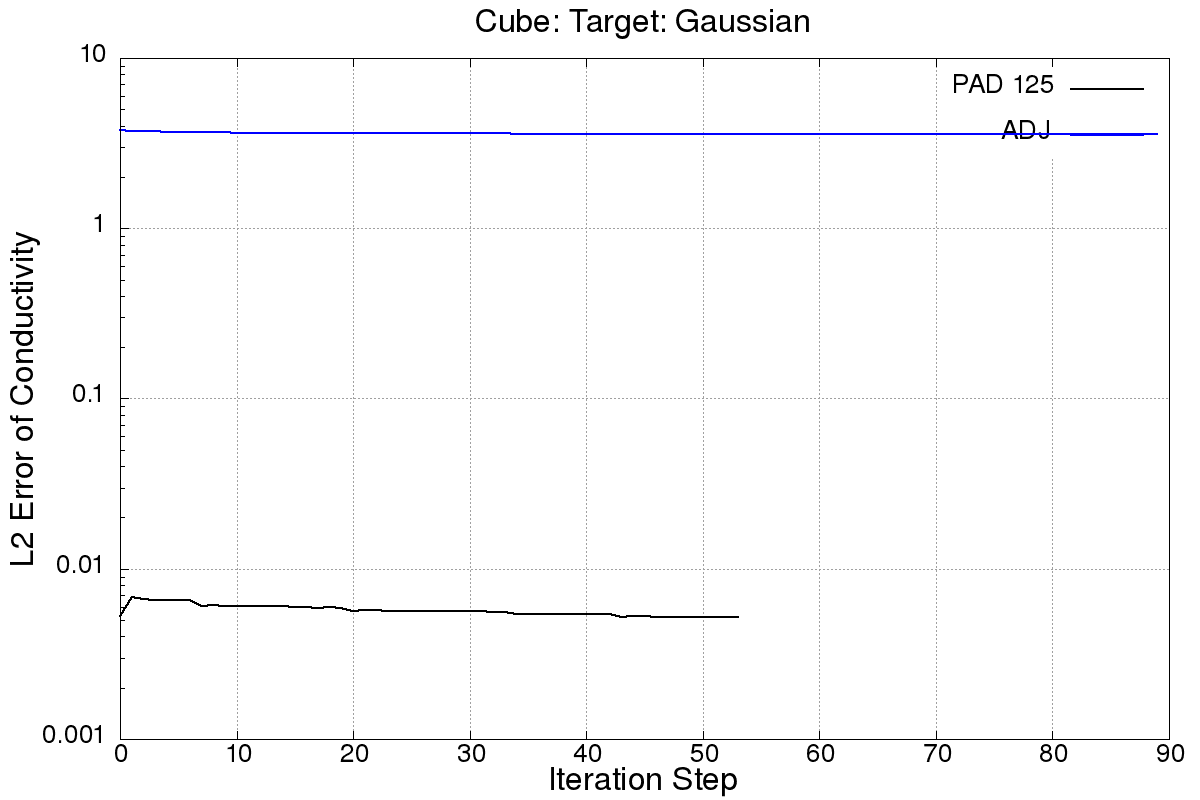}	
	\caption{\label{f:cubeGauss8} Cube - Gaussian $k$: Errors in Fluxes (left) and Conductivity (right).}
\end{figure}

\noi
Three measurement functions similar to the ones used before were
placed at the center of each face so as to induce a heat flux that
is mainly directed in the $x,y,z$-directions repectively.
For the optimization via finite differences, the distribution
of $k$ was given by constant regions spaced in a cartesian lattice in
$x,y,z$. Figures~\ref{f:cubeGauss}-\ref{f:cubeGauss6} show the results 
obtained using 125 degrees of freedom for the finite difference gradient and 
the more than 250,000 degrees of freedom when using each element
conductivity. Figure~\ref{f:cubeGauss7} depicts the comparison of fluxes for 
plane $y=0$. One can see that the optimization procedure as
such works, recovering well the specified normal flux distribution on
the boundary, even though the spatial distribution for $k$ is very
different from the target. As before, apparently the lower the 
degrees of freedom, the better the result (!). The decrease in 
objective function for these cases is shown in 
Figure~\ref{f:cubeGauss8} (left), and the L2 error of the conductivity 
in Figure~\ref{f:cubeGauss8} (right). 

\clearpage

\section{Conclusions and Outlook}

A finite element code for heat conduction, together with an
adjoint solver and a suite of optimization tools were applied
for the solution of Calderon's problem. One of the questions
whose answer was sought was whether the solution to these
problems is unique and obtainable. The results to date show
that while the optimization procedure is able to obtain spatial
distributions of the conductivity $k$ that reduce the cost function
significantly, the resulting conductivity $k$ is still significantly
different from the target distribution sought, particularly if a
single temperature and flux field (i.e. `measurement') is considered.
While the normal fluxes recovered are very close to the prescribed 
ones, the tangential fluxes can differ considerably.
As a possible way to circumvent these difficulties the prescription 
of several temperature and flux fields (i.e. `measurements') at 
the same time were explored. The aim was to have as 
many boundary regions as possible with significant normal fluxes. 
The observation made here is that this technique yielded closer
spatial distributions of the conductivity $k$ to the targets
desired, but that the effectiveness of number of `measurements' 
tends to saturate at a relatively low number. For the cases run
to date 2-4 measurements in 2-D and 3-6 measurements in 3-D seemed
to give the best results.

At this point, it is not clear why rigorous mathematical proofs
yield results of convergence and uniqueness, while in practice
accurate distributions of the conductivity $k$ seem to be elusive.
One possible explanation is that the spatial influence of 
conductivities decreases exponentially with distance. Thus, many 
different conductivities inside a domain could give rise to very similar
(infinitely close) boundary measurements.

\bibliographystyle{plain}
\bibliography{refs}

\end{document}